\documentclass[10pt]{amsart}

\usepackage{amssymb} 
\usepackage[latin1]{inputenc}
\usepackage[matrix,arrow,curve]{xy}
\usepackage[mathscr]{eucal}
\usepackage{epic,eepic}
\usepackage{bbm,url}

\makeindex
\newtheorem{theorem}{Theorem}[section]
\newtheorem{itheorem}{Theorem}
\newtheorem{icorollary}{Corollary}

\newtheorem{lemma}[theorem]{Lemma}
\newtheorem{prop}[theorem]{Proposition}

\newtheorem{cor}[theorem]{Corollary}

\theoremstyle{definition}
\newtheorem{idefinition}{Definition}
\newtheorem{iremark}{Remark}
\newtheorem{definition}[theorem]{Definition}
\newtheorem{remark}[theorem]{Remark}
\newtheorem{example}[theorem]{Example}

\newenvironment{pf}
{\medskip\noindent {\it Proof --- \ }}
{\hfill\nobreak $\Box$ \par\bigbreak}

\renewcommand{\L}{{\mathcal L}}

\newcommand{\RR}{{\mathcal R}}

\newcommand{\Hom}{\text{Hom}}

\newcommand{\F}{ \mathbb F} 
\newcommand{\C}{{ \mathbb C  }}

\newcommand{\Q}{{ \mathbb Q } }
\newcommand{\Z}{{ \mathbb Z  }}
\newcommand{\N}{{ \mathbb N  }}

\newcommand{\HH}{{\mathcal H}}
\newcommand{\PP}{{\mathbb P}}

\newcommand{\Gal}{{\mathrm{Gal}\,}}

\newcommand{\A}{{\bf A}}

\newcommand{\End}{{\text{End}}}
\newcommand{\Id}{{\text{Id}}}

\newcommand{\U}{{\text{U}}}

\newcommand{\Gl}{{\text {GL}}}

\newcommand{\GL}{{\text {GL}}}

\newcommand{\Sl}{{\text {SL}}}

\newcommand{\T}{{\mathbb T}}
 
\newcommand{\spec}{{\text{Spec\,}}}

\newcommand{\TT}{{\mathbb T}}

\newcommand{\vareps}{{\varepsilon}}
\newcommand{\CC}{{\mathcal{C}}}

\newcommand{\D}{{\bf D}}

\newcommand{\G}{{\mathfrak g}}

\newcommand{\ad}{{\text{ad}}}

\newcommand{\Ext}{{\text{Ext}}}

\newcommand{\un}{{1\hspace{-0.15cm}1}}

\newcommand{\ses}{{\text{ss}}}

\newcommand{\diag}{{\text{diag}}}
\newcommand{\Frob}{{\text{Frob\,}}}

\newcommand{\tr}{{\text{tr\,}}}

\newcommand{\rhob}{{\bar \rho}}

\renewcommand{\hat}{\widehat}

\newcommand{\isomo}{\overset{\sim}{\longrightarrow}}

\renewcommand{\v}{{\bf v}}

\newcommand{\W}{W_{\Q_p}}

\renewcommand{\Gal}{{\rm Gal}}

\newcommand{\cont}{\mathrm{cont}}

\renewcommand{\G}{G}
\renewcommand{\sp}{{\mathrm{Sp}}}
\newcommand{\crys}{\text{crys}}

\newcommand{\WW}{{\mathcal{W}}}
\newcommand{\hotimes}{\hat \otimes}
\newcommand{\teich}[1]{\langle#1\rangle}
\newcommand{\Symb}{{\rm{Symb}}}

\newcommand{\tors}{\text{tors}}
\newcommand{\tw}{{\tilde{w}}}

\renewcommand{\diamond}{\text{diamond operators }}
\newcommand{\diamonds}{\diamond $\langle a \rangle$ ($a \in (\Z/N\Z)^\ast$)}
\newcommand{\free}{\text{free}}
\newcommand{\mb}{\par \medskip}
\newcommand{\es}{\text{es}}
\newcommand{\Mel}{\text{Mel}}
\newcommand{\NN}[1]{|\!| #1 |\!|}
\newcommand{\spe}{\text{sp}}
\begin{document}

\baselineskip 15.8pt

\bibliographystyle{style} 
\title{Critical $p$-adic $L$-functions}
\date{November 2, 2009}

\begin{abstract} 
We attach $p$-adic $L$-functions to critical modular forms and study them. We prove that those $L$-functions fit in a two-variables
$p$-adic $L$-function defined locally everywhere on the eigencurve. 
\end{abstract}

\thanks{During the elaboration of this paper, I was supported by the NSF grant DMS 08-01205.
I thank Ga\"etan Chenevier for his invaluable help related to this article. I can say that most of section 2 and 3 below
arose from discussions with him: not only did he give me many individual ideas or arguments, but I owe him my whole vision of the subject.
I also thank Glenn Stevens and Robert Pollack for giving me access to their work even before it was made available on the web,  
and for many conversations about the themes of this paper. Also the remarkable clarity of their arguments in the various papers by them in the bibliography were a strong appeal for me to study this subject, and a  great help when I did. Finally, I thank Kevin Buzzard, Matthew Emerton and Lo\"ic Merel for useful conversations.}

\author[J.~Bella\"iche]{Jo\"el Bella\"iche}
\email{jbellaic@brandeis.edu}
\address{Jo\"el Bella\"iche\\Brandeis University\\
415 South Street\\Waltham, MA 02454-9110\\U.S.A}

\maketitle

%\tableofcontents

\tableofcontents

\section{Introduction}

\subsection{Objectives and motivation}

The aim of this paper is to extend, and even in a sense, to complete the works of Mazur-Swinnerton-Dyer (\cite{msd}), 
Manin (\cite{manin}), Visik (\cite{visik}), Amice-V\'elu (\cite{amicevelu}), Mazur-Tate-Teiltelbaum (\cite{MTT}), Stevens (\cite{stevens}), Stevens-Pollack (\cite{stevenspollack1} and \cite{stevenspollack2}) on the construction of a $p$-adic $L$-function $L(f,s)$ of one $p$-adic 
variable $s$ attached to a refined classical modular form $f$ that is an eigenvector for almost all Hecke operators,
as well as the work of Mazur, Kitagawa (\cite{kitagawa}), Greenberg-Stevens (\cite{greenbergstevens}), Stevens
\cite{stevensfam}, Panchishkin (\cite{panchishkin}),  and Emerton (\cite{emerton}) on the construction of a two-variables $p$-adic $L$-function $L(x,s)$ where the first variable $x$ runs among the points in a suitable 
$p$-adic family of refined modular forms, such that
when $x$ corresponds to a classical refined modular eigenform $f$, the function $s \rightarrow L(x,s)$ is the classical one-variable $p$-adic $L$-function $L_f(s)$. 

\par \medskip
Let us give a little bit more details. It has been known for more than three decades
how to attach $p$-adic $L$-functions to refined modular forms of non-critical slope. More recently, Stevens and Pollack
have constructed $L$-functions for certain refined modular form of critical slope, namely the ones that are not $\theta$-critical.
Their construction therefore does not work for critical slope CM cuspidal forms and Eisenstein series, which are very important for number-theoretic applications. In this paper, we construct in a unified way a $p$-adic $L$-function for (almost) all refined modular forms.
Our method uses the insights of Stevens and Pollack-Stevens ($p$-adic $L$-function as coming from distributions-valued
modular symbols), together with a knowledge of the geometry of the eigencurve and of various modules over it
in the spirit of \cite{BCsmooth}. So we use family arguments even for the construction of the $L$-function of an individual
refined modular form. 

Our second main result is the construction of a two-variables $L$-function $L(x,s)$ interpolating the individual $L$-functions
of refined modular form (locally on $x$ in the neighborhood of (almost) any refined modular form in the eigencurve). The existence of such a function restricted to the $x$'s of critical slope have been shown by three
different authors (firstly Stevens in a work \cite{stevensfam} that unfortunately was not made available, then Panchishkin \cite{panchishkin}  and Emerton \cite{emerton}).
We extend their results by incorporating the missing $x$'s. Our method is very close to what I believe is Stevens' method, using families
of distributions-valued modular symbols over the weight space,
but there is  a new difficulty, which is that the eigencurve is not \'etale over the weight space at the critical point $x$.  In the cases 
already treated in the literature, 
the \'etaleness of the eigencurve at the considered point made that, locally, the variable $x$ on the eigencurve was essentially 
the same as the variable $k$ on the weight space (and two-variables $L$-functions are accordingly written $L(k,s)$
rather than $L(x,s)$ in the existing literature.)  In our case, we have to use our knowledge of the geometry of the eigencurve again to overturn this difficulty.

The interplay between the geometry of the eigencurve (absolute and relative to the weight space)
and the construction and properties of the $L$-functions is promising of interesting phenomena in higher rank, when the
geometry of the eigenvarieties can be much more complicated (see e.g. \cite{Bduke}) and of number-theoretic significance 
(see \cite[Chapter IX]{BCbook}). For more on this, see~\ref{perspectives} below.

\par \bigskip
In the remaining of this introduction, I will explain in more details my results and their relations 
with the results obtained by the authors  aforementioned, especially with  the results of Stevens and Pollack, 
 which were a fundamental source of inspiration for me. For this I need to introduce some notations and terminology.

\subsection{Notations and conventions}
\label{notations}

Throughout  this paper, we shall fix an integer $N \geq 1$, and a prime number $p$ that we shall assume\footnote{Let us just mention, for the sake of completeness, what happens if $p$ divides $N$. Then a newform $f$ for $\Gamma_1(N)$ has only one refinement, itself.
The  
eigenvalue $\alpha$ of $f$ for the operator $U_p$ satisfies either $v_p(\alpha)=(k+1)/2$ or $\alpha=0$. In the first case,
$f$ is of non-critical slope and the construction of its $p$-adic $L$-function was known back in the seventies. 
In the second case, the slope is infinite, and the situation is very mysterious. Our methods, relying on seeing $f$ 
as a point on the eigencurve, does not apply and we have nothing new to say in this case.}
 odd and prime to $N$.
In all normed extension of $\Q_p$ we shall use the normalized valuation $v_p$ ($v_p(p)=1$) and the normalized absolute value $p^{-v_p}$.
We fix an embedding $\bar \Q \hookrightarrow \bar \Q_p$ and $\bar \Q \hookrightarrow \C$.

We shall work with the congruence subgroup of $\Sl_2(\Z)$ called $\Gamma$ defined by $\Gamma=\Gamma_1(N) \cap \Gamma_0(p)$.
For an integer $k \geq 0$, we shall denote by $M_{k+2}(\Gamma)$ (resp. $S_{k+2}(\Gamma)$, $M^\dag_{k+2}(\Gamma)$, $S^\dag_{k+2}(\Gamma)$) the $\Q_p$-spaces of classical modular forms 
(resp. cuspidal classical modular forms, resp. overconvergent $p$-adic modular forms, resp. overconvergent cuspidal $p$-adic modular form) of level $\Gamma$ and weight $k+2$. Those spaces is acted upon by the Hecke operators $T_l$ for 
$l$ prime to $Np$, the Atkin-Lehner operator $U_p$, and the diamond operators $\langle a \rangle$ for $a \in (\Z/N\Z)^\ast$, as will be 
virtually all spaces and modules considered in this paper. To avoid repeating this list of operators too often, we shall define $\HH$ as the commutative polynomial algebra over $\Z$ over the variables $T_l$ (for $l$ prime to $Np$), $U_p$, and $\langle a \rangle$ for $a \in (\Z/\N\Z)^\ast$, that is
$$\HH=\Z[(T_l)_{l \not\, | pN}, U_p, ( \langle a \rangle)_{a \in (\Z/\N\Z)^\ast}],$$
and simply say
that the spaces $M_{k+2}(\Gamma)$ and $S_{k+2}(\Gamma)$ are acted upon by $\HH$:

In general, if $M$ is a $\Q_p$-space on which $\HH$ acts, by a {\it system of $\HH$-eigenvalues} appearing in $M$ we shall mean
any character $x: \HH \rightarrow \bar \Q_p$ such that there exists a vector $v$ in $M \otimes \bar \Q_p$ satisfying $h v = x(h) v$ for all $h \in \HH$. We call such a $v$ an $x$-eigenvector. If $x$ is such a system, we shall  call the {\it eigenspace for the system of $\HH$-eigenvalues $x$} or shortly the {\it $x$-eigenspace} the space, denoted by $M[x]$, of $x$-eigenvectors $v$ in $M \otimes \bar \Q_p$. We shall call  {\it generalized eigenspace for the system of $\HH$-eigenvalues $x$}, denoted by $M_{(x)}$, or shortly the generalized $x$-eigenspace
 the space of $v$ in $M \otimes \bar \Q_p$ such that for all $h \in \HH$, there exists an $n \in \N$ such that $(h-x(h))^n v=0$. Please keep in mind that we have extended scalar to $\bar \Q_p$, so that $M[x]$ and $M_{(x)}$ are 
$\bar \Q_p$-vector spaces. If $N$ is another (or the same) $\Q_p$-spaces of which $\HH$ acts, and $v \in N$ is an $\HH$-eigenvector,
of system of eigenvalues $x$, we shall write $M[v]$ and $M_{(v)}$ instead of $M[x]$ and $M_{(x)}$.

\subsection{Reminder about modular forms and the eigencurve}

To state our results, and to explain in more details their relations with the works mentioned above,
we need to recall the notions of a non necessarily cuspidal newform, of  a {\it refinement} of a modular form, of a {\it critical slope},
$\theta$-critical, and {\it critical} refinement.

 Let $f=\sum a_n q^n$ be a modular form
of level $\Gamma_1(N)$ and of weight $k+2 \geq 2$ with coefficients in $\bar \Q_p$.
We shall say that $f$ is a {\it newform} (of weight $k+2$, level $\Gamma_1(N)$ ) if $f$ is eigenvector for the Hecke operators $T_l$ ($l$ prime to $N$) and the \diamonds\ , and if there is no modular form of weight $k+2$ and level $\Gamma_1(N')$ where $N'$  is a {\bf proper divisor of $N$} which is also an 
eigenform for all the $T_l$, $(l,N)=1$, and all the \diamond\  {\bf with the same eigenvalues as $f$.}

If $f$ is newform of weight $k+2$, level $\Gamma_1(N)$, we have $T_p f = a_p f$ for some algebraic integer $a_p$ and 
$\langle p \rangle f = \epsilon(p) f$ for some root of unity $\epsilon(p)$. Let us call $\alpha$ and $\beta$ the roots of the polynomial
\begin{eqnarray} \label{hp} X^2-a_p X + \epsilon(p) p^{k+1}.\end{eqnarray}
We shall assume throughout the paper without loss of generality\footnote{This is conjectured in general, known in weight $k+2=2$, and when this is not true we have $v_p(\alpha)=v_p(\beta)=(k+1)/2 < k+1$, so the results we want to prove are already known.} 
 and for simplicity that $\alpha \neq \beta$.
As is well known, there are two  normalized modular forms for $\Gamma=\Gamma_1(N) \cap \Gamma_0(p)$, of level $k+2$, that are eigenforms for the Hecke operators $T_l$ ($l$ prime to $pN$) and the diamond operators with the same eigenvalues as $f$, and that are also eigenforms 
for the Atkin-Lehner operator $U_p$:
\begin{eqnarray}  f_\alpha(z) &=& f(z)-\beta f(pz) \\
\label{fbeta} f_\beta(z) &=& f(z) - \alpha f(pz).
\end{eqnarray}
Those forms 
satisfy $U_p f_\beta = \beta f_\beta$ and $U_p f_\alpha = \alpha f_\alpha$.

The forms $f_\beta$ and $f_\alpha$ are called {\it refined modular forms}, and they are the two 
{\it refinements of $f$}. The choice of a root $\beta$ or $\alpha$ of (\ref{hp}) is also called a {\it refinement of $f$}. Refinements are also sometimes
called {\it $p$-stabilizations}, but we shall not use this terminology. 
Refined modular form are the natural objects to which attach a $p$-adic
$L$-function, and they are also (a coincidence due to our two-dimensional setting) the natural object to put in $p$-adic families:
in particular a refined modular form $f_\beta$ corresponds to a point $x$ in the $p$-adic eigencurve $\CC$
of tame level $\Gamma_1(N)$. We call such points the {\it classical points of $\CC$}.
Hence classical points of $\CC$ are in bijection with refinements of modular newforms of level $\Gamma_1(N)$ and arbitrary weight $k+2$. 

If $f_\beta$ is a refined modular form as above, we shall say that it is of {\it critical slope} if $v_p(\beta)=k+1$, is {\it $\theta$-critical}
if $f_\beta$ is in the image of the operator $\theta^{k+1}$, and is {\it critical} if $M_{k+2}(\Gamma)_{(f_\beta)} \neq M_{k+2}^\dag(\Gamma)_{(f_\beta)}$.    We refer the reader to ~\S\ref{refinement} for equivalent conditions and more details. Let us just say now for the intelligence of our results that critical implies $\theta$-critical implies critical slope; that for a cuspidal form, critical and $\theta$-critical are equivalent; and that for an Eisenstein series,
$\theta$-critical and critical slope are equivalent; otherwise implications are strict.

\subsection{Results}

\subsubsection{Modular symbols}

\label{intmodsymb}

To formulate our first result about the existence of a $p$-adic $L$-function, we shall adopt Stevens' point of view
of $L$-functions as Mellin transforms of the values at the divisor $\{\infty\}-\{0\}$ of distributions-valued modular symbols. 
We begin by recalling quickly this point of view.

Let $0<r\leq 1$ be any real number, and let $\D=\D[r]$ be the $\Q_p$-Banach space which is the dual of the Banach space of functions on $\Z_p$ that are analytic on each closed ball of radius $r$. It will be
convenient to take $r=1/p$ in this introduction, but actually nothing would change for any other value of $r$. 
 
For a fixed integer $k$, $\D$ is endowed with an action of $\Gamma_0(p)$, called the {\it weight-$k$ action}, which is the dual of the action on the space of locally analytic functions given for $\gamma = \left( \begin{matrix} a & b \\ c & d \end{matrix}\right)$
by $$(\gamma \cdot_k f)(z)  =  (a+cz)^k f\left(\frac{b+dz}{a+cz}\right).$$ For $\Gamma=\Gamma_1(N) \cap \Gamma_0(p)$,
we define 
the space of modular symbols $\Symb_\Gamma(\D_k):=\Hom_\Gamma(\Delta_0,\D_k)$ where $\Delta_0$ is the abelian
group of divisors of degree $0$ on $\PP^1(\Q)$. This $\Q_p$-Banach space
is endowed with a natural $\HH$-action and
also with an involution $\iota$ (commuting with $\HH$) given by the matrix $\left( \begin{matrix} 1 & 0 \\ 0 & -1 \end{matrix} \right)$. We denote by  $\Symb_\Gamma^{\pm}(\D_k)$ the eigenspaces of this involution. For more details about those notions see below \S\ref{remindermodsymb}. Also, $\D_k$ admits as a quotient the space $V_k$ of 
homogeneous polynomial of degree $k$ in two variables, and there is an $\HH$-equivariant "specialization" morphism 
$\rho_k^\ast: \Symb_\Gamma^\pm(\D_k) \rightarrow \Symb_\Gamma^\pm(V_k)$. 
For more details about those notions see below \S\ref{stevenscontrol}.

Our theorem has a technical condition, but fortunately it is very mild.
\begin{idefinition} \label{defdecent} We shall say that a refined form $f_\beta$ is {\it decent} if it satisfies {\bf at least one} of the three following assumption:
\begin{itemize}
\item[(i)] $f$ is Eisenstein
\item[(ii)] $f_\beta$ is non-critical.
\item[(iii)] $f$ is cuspidal and $H^1_g(G_\Q,\rho_f)=0$
\end{itemize} 
\end{idefinition}
As we shall see in~\S\ref{decent}, all CM forms are decent, and there exist several general easy to check  criteria ensuring that a cuspidal non-CM form
satisfies (iii) (so is decent) leaving only a small number of residual forms out. 
Moreover, it is conjectured that every cuspidal form satisfies (iii), and independently that every $f$ which is cuspidal non CM satisfies (ii), so
cuspidal non-CM forms should be doubly decent.
 
\begin{itheorem} \label{eigenspace}
Let $f$ be a newform of weight $k+2$ and level $\Gamma_1(N)$, and let
$f_\beta$ a refinement of $f$ at $p$. We assume that $f_\beta$ is decent, and if $f$ is Eisenstein, that $v_p(\beta) > 0$.
Then for both choices of the sign $\pm$, the  eigenspace $\Symb_\Gamma^\pm(\D_k)[f_\beta]$ has dimension $1$.

Moreover the following are equivalent 
\begin{itemize} \item[(i)] $f_\beta$ is critical.
\item[(ii)] The generalized eigenspaces  $\Symb_\Gamma^\pm(\D_k)_{(f_\beta)}$ (which always have the same dimension for both choices of the sign $\pm$) have dimensions strictly greater than $1$. 
\item[(iii)] One has
$$ \text{(TS${}_{\pm}$)\ \ \ \ \ \ \ \ \ } \rho_k^\ast(\Symb_\Gamma^\pm(\D_k)[f_\beta]) = 0.$$ for {\bf both} choice of the sign $\pm$
\end{itemize} 
If $f$ is cuspidal, then if (TS${}_{\pm}$) holds for one choice of the sign $\pm$
it holds for the other as well. If $f$ is Eisenstein, (TS${}_{-\epsilon(f_\beta)}$) 
always holds if $\epsilon(f_\beta)$ is the sign of $f_\beta$ defined in \S\ref{signeisenstein}.
\end{itheorem} 
This theorem was known by a work of Stevens-Pollack (\cite{stevenspollack2}) in the case of forms $f_\beta$ that are not
$\theta$-critical and before in the case of  forms of non-critical slope by a result of Stevens (\cite{stevens}). Therefore, the new cases provided of this theorem are the $\theta$-critical cases, which include 
 all Eisenstein series (with their non-ordinary refinement), all CM forms (such that $p$ is split in the quadratic field attached to $f$,
and with their non-ordinary refinement), and the cuspidal non CM $f_\beta$ that are
not known to be non-$\theta$-critical.
 
 \begin{iremark}
 We note that the the one-dimensionality of the eigenspace $\Symb_\Gamma^\pm(\D_k)[f_\beta]$ was observed in several instances using numerical computations with a computer by Pasol, Pollack, and Stevens. For example, for the modular form attached to the CM elliptic curve $X_0(32)$ with $p=5$ (\cite[Example 6.10]{stevenspollack2}, \cite{stevenspollack1}), or for the Eisenstein series of level $\Gamma_0(11)$ and weight $2$ for $p=3$ (\cite{stevensEisenstein}). It was furthermore observed that the unique eigensymbol corresponding to the above CM form specializes (by the morphism $\rho_k^\ast$, see above) to $0$, while this was not the case for the above Eisenstein series. It was even conjectured (\cite["Wild Guess"]{stevensEisenstein})  that the same result would hold for Eisenstein Series for $\Gamma_0(11)$ of all weight, and $p=3$. Those observations and conjectures
 are consequences of our result since critical slope CM forms are always critical, while critical slope Eisenstein series are never critical for regular primes $p$ -- see Prop~\ref{criticalclass}.
 
 Those observations and conjectures were important motivations for the present work.
 \end{iremark}

\subsubsection{Construction of $p$-adic $L$-functions}

The connection with the $p$-adic $L$-functions is as follows: 
Let $\Phi_{f_\beta}^\pm$ be generators of the eigenspaces  $\Symb_\Gamma^{\pm}(\D_k)[f_\beta]$
considered in Theorem~\ref{eigenspace}. Both generators are well defined up to multiplication by a non-zero scalar (in $\bar \Q_p^\ast$).
We define two $p$-adic $L$-function  $L^\pm(f_\beta,\sigma)$ of the refined form $f_\beta$  as analytic functions 
of a continuous character $\sigma : \Z_p^\ast \rightarrow \C_p^\ast$ by 
$$L^\pm (f_\beta,\sigma) := \Phi_{f_\beta}^\pm(\{\infty\}-\{0\})(\sigma).$$
Here by definition $\Phi^\pm_{f_\beta}(\{\infty\}-\{0\})$ is a linear form on the space of locally analytic function on $\Z_p$ of radius of convergence $r=1/p$ everywhere and  $\Phi_{f_\beta}^\pm(\{\infty\}-\{0\})(\sigma)$ 
is the value of that linear form on $\sigma$, which is seen as the locally analytic function on $\Z_p$ which is the character $\sigma$ on $\Z_p^\ast$ and $0$ on $p\Z_p$. The two functions $L^\pm(f_\beta,\sigma)$ are each defined up to multiplication by a non-zero scalar, but note that 
$L^\pm (f_\beta,\sigma)=0$ whenever $\sigma(-1) \neq \pm 1$, so the two functions have disjoint supports.
Therefore the zeros and poles of the following function
$$L(f_\beta,\sigma):=L^+(f_\beta,\sigma) + L^-(f_\beta,\sigma)$$
are well-defined, even if $L(f_\beta,\sigma)$ depends on two auxiliary scalars.

When $f_\beta$ is non-critical, this $p$-adic $L$-function is the same as the usual one, defined by Mazur-Tate-Teitelbaum in
the non-critical slope case, and by Stevens-Pollack in the critical slope non-$\theta$-critical case. This is clear by definition,
at least if we use Steven's reformulation (\cite{stevens}) of the definition of the $L$-function in the non-critical slope case. In other words, our construction is new only in the $\theta$-critical case.

\subsubsection{Properties of the $p$-adic $L$-functions}
\label{properties}

The $p$-adic $L$-functions enjoy the following properties: 
\mb

{\bf Analyticity:} 
The function $L(f_\beta,\sigma)$ is an analytic function of $\sigma$, in the following sense:

Let $\mu_{p-1}$ be the group of $(p-1)$-th roots of unity in $\Z_p^\ast$. As is well known,
$\mu_{p-1}$ is cyclic of order $p-1$, isomorphic through the reduction mod $p$ to $\F_p^\ast$, and there is a canonical isomorphism $\Z_p^\ast \isomo \mu_{p-1} \times (1+p \Z_p)$ sending $x$ to $(\teich{x},x/\teich{x})$ where $\teich{x}$ is the Teichmuller representative of $x \mod p$ in $\mu_{p-1}$. 

We fix a generator $\gamma$ of $1+p \Z_p$ and for $s \in \C_p,\  |s-1|<1$, let $\chi_s$ be the continuous character of $1+p \Z_p$ sending $\gamma$ to $s$. All continuous $\C_p$-valued characters of $1+p \Z_p$ are of this form. 

Therefore, any continuous character $\sigma : \Z_p^\ast \rightarrow \C_p^\ast$ can be written in a unique way $\chi=\psi \chi_s$  where $\psi$ is a character of $\mu_{p-1}$
and $\chi_s$ a character of $1+p \Z_p$, with $|s-1|<1$.

The assertion that $L(f_\beta,\sigma)$ is analytic in $\sigma$ is to be understood as follows: For every character $\psi$ of $\mu_{p-1}$, 
the function $s \mapsto L(f_\beta,\psi \chi_s)$ is given
by a power series in $s$ converging on the open ball $|s-1|<1$.

\mb

{\bf Order:} The function $L(f_\beta,\sigma)$ is of order at most $v_p(\beta)$, in the sense that
we have $L^\pm(f_\beta,\psi \chi_s)=O( (\log_p(s))^{v_p(\beta)})$, where for $F$ and $G$ two functions of $s$ defined on the set $|s-1|<1$, we  have $\sup_{|s-1|<\rho} |F(s)| = O(\sup_{|s-1|<\rho} G(s))$ when $\rho \rightarrow 1^{-}$.

\mb
{\bf Interpolation:} We know the value of $L(f_\beta,\sigma)$ at all characters $\sigma$ of the form $t \mapsto \varphi(t) t^j$, where $\varphi$ is a character of finite order
 of $\Z_p^\ast$,  and $j$ is a rational integer such that $0 \leq j \leq k$.
When $v_p(\beta)<k+1$, those values determine uniquely the function
$L^\pm(f_\beta,\sigma)$ in view of the above property. This is not true however in the critical-slope case $v_p(\beta)=k+1$.

If $f_\beta$ is not $\theta$-critical, the interpolation property is due to Pollack and Stevens (or earlier works if $v_p(\beta)<k+1$):
\begin{eqnarray} \label{intnoncrit} L^\pm(f_\beta,\varphi t^j)=e_p(\beta,\varphi t^j) \frac{m^{j+1}}{(-2\pi i)^j} \frac{j!}{\tau(\chi^{-1})} \frac{L_\infty(f\varphi^{-1},j+1)}{\Omega_f^{\pm}}\end{eqnarray}
 for all $\varphi : \Z_p^\ast \rightarrow \Q_p^\ast$ of conductor $m=p^\nu M$ (with $M$ prime to $p$) satisfying $\varphi(-1)=\pm 1$ and $0 \leq j \leq k$, up to a scalar independent on $\varphi t^j$. Above, $L_\infty$ is the archimedean $L$-function, $\Omega^{\pm}$ are the two archimedean periods
of $f$, $\tau(\varphi^{-1})$ is the Gauss sum of $\varphi^{-1}$ and 
\begin{eqnarray} \label{ep}
e_p(\beta, \varphi t^j)= \frac{1}{\beta^\nu}  (1 - \frac{\bar \varphi(p) \epsilon(p) p^{k-j}}{\alpha}) (1 - \frac{\varphi(p) p^j}{\alpha})
\end{eqnarray}

\begin{itheorem}\label{thminterpol}
If $f_\beta$ is $\theta$-critical, then
\begin{eqnarray} \label{intcrit}
L^\pm(f_\beta,\varphi t^j) = 0 \end{eqnarray} for all $\varphi: \Z_p^\ast \rightarrow \Q_p^\ast$
of finite order and $0 \leq j \leq k$, excepted in the following 
case: $f$ is Eisenstein, $f_\beta$ is not critical\footnote{This condition is expected to always hold for $f$ Eisenstein, see Prop.~\ref{criticalclass}.}, 
the sign $\pm$ is $\epsilon(f_\beta)$, the character $\chi$ is trivial, and $j=k$.
% In the exceptional case, one has $L^{\epsilon}(f_\beta,t^k) \neq 0$.
\end{itheorem} 
One can reformulate this theorem by introducing the analytic function of the variable $s$
(on the open unit ball of center $1$ and radius $1$, that is $|s-1|<1$)
 $$\log^{[k]}(s) := \prod_{i=0}^k \log_p (\gamma^{-j} s),$$
 whose list of zeros, all with simple order, are all the points $s = \gamma_j \zeta$, where $\zeta$ is a $p^\infty$-root of unity and $0 \leq j \leq k$. Then the theorem simply says
that {\it if $f_\beta$ is $\theta$-critical, then 
for all tame characters $\psi$, the analytic function $s \mapsto 
L^\pm(f_\beta,\psi \chi_s)$ is divisible by $\log^{[k]}(s)$ in the ring of analytic function 
on the open unit ball of center $1$, excepted in the exceptional
case where $f$ is Eisenstein, $f_\beta$ non-critical, $\psi:\mu_{p-1} \rightarrow \C_p^\ast$ is the raising to the $k$-th power  $\psi_k(t) = t^k$, $\pm=\epsilon(f_\beta)$; in that case the function 
$L^{\epsilon(f_\beta)}(f_\beta, \psi_k \chi_s)$ is divisible by $\frac{\log^{[k]}(s)}{s-\gamma^k}$.}
(To see that this is equivalent to the theorem, it is enough to observe that a special character 
$t \mapsto \varphi(t) t^j$, when written in the canonical form $\psi \chi_s$, has 
$\psi(t) = \varphi(t) t^j$ for $t \in \mu_{p-1}$ and $s=\varphi(\gamma) \gamma^j$, where $\varphi(\gamma)$ is a $p^\infty$-root of unity. In particular the character $ \psi_k \chi_{\gamma^k} $ is simply the character $t \mapsto t^k$.)

As a consequence of this and the {\bf order} property, 
if $f_\beta$ is $\theta$-critical, then $L^\pm(f_\beta,\psi \chi_s)/\log^{[k]}(s)$ is a meromorphic function on the open ball $|s-1|<1$ with a finite number of $0$'s
and at most one simple pole at $s=1$. 

\par \bigskip
Let us recall that it is conjectured that $L(f_\alpha,\psi \chi_s)$ has infinitely many $0$'s whenever $\v_p(\alpha)>0$. This is known (cf. \cite{pollack}) in the case where the coefficient $a_p$ of the modular form $f$ is $0$ (e.g. a modular form attached to a super-singular elliptic curve).
\begin{icorollary} \label{corinterpol}
If $v_p(\beta)=k+1$ (that is $f_\beta$ of critical slope), then for all character $\psi : \mu_{p-1} \rightarrow \C_p^\ast$, the function $s \mapsto L^\pm(f_\beta,\psi \chi_s)$
has infinitely many zeros.
\end{icorollary}

\begin{iremark}
\begin{itemize}
\item[(i)] Since all the new (that is, $\theta$-critical) $p$-adic $L$-functions
$L^\pm(f_\beta,\psi\chi_s)$ are divisible by $\log^{[k](s)}$, it is  natural to try to compute the quotients $$L^\pm(f,\beta)(\psi\chi_s)/\log^{[k]}(s).$$ Let us focus on the cases where $f$ is CM or Eisenstein (conjecturally, that should be the only cases 
where there is a $\theta$-critical $f_\beta$). It is natural to expect for 
 $L^\pm(f,\beta)(\psi\chi_s)/\log^{[k]}(s)$ a formula in terms of $p$-adic $L$-functions of 
 Hecke characters of $\Q$ (that is Kubota-Leopoldt's $L$-functions) in the case of an Eisenstein series, and of the quadratic imaginary field $K$ (that is Katz' $L$-functions)
 in the case of a form CM by $K$. 
 
 Actually, in three cases out of four, such a formula is elementary to prove.
 Those cases are the cases where the modular symbols $\Phi_{f_\beta}^\pm$
 is in the kernel of $\rho_k^\ast$ (hence in the image of the mat $\Theta_k^\ast$ defined in \S\ref{stevenscontrol} below), namely the two cases corresponding to both values of the sign $\pm$ if $f$ is CM, and the case $\pm = - \epsilon(f_\beta)$ if $f$ is Eisenstein. 
 
 The remaining case, namely the case of $L^{\epsilon(f_\beta)}(f_\beta,\psi\chi_s)$ is much more mysterious. In the special case $p=3$, $k=2$, $N=11$, it is the object of a conjecture of Stevens (\cite{stevensEisenstein}) already mentioned above. In a work in progress with several other mathematicians, I plan to solve this case as well. This work shall also contain the details of the computations in the { easy cases}.
 
 A closely related question is to prove a suitably formulated Iwasawa conjecture relating the "analytic" $L$-function
 $L(f_\beta,\psi \chi_s)$ and an "algebraic" L-function attached to the preferred refined
 representation $\rho_{f_\beta}$ (see \S\ref{preferred} below).
 
 \item[(ii)] Do we know that our new $L$-functions $L^\pm(f_\beta,\cdot)$ are not identically zero?
 Before answering this question, I want to stress that it is  fundamentally different 
 from asking whether the distributions-valued modular symbol $(\Phi_{f_\beta}^\pm)$  
 defining the $L$-function is $0$ (it is not, by definition and Theorem~\ref{eigenspace}).
 Indeed, we see this difference even in the classical (that is not $\theta$-critical) case, 
 since there the non-vanishing of $\Phi_{f_\beta}^\pm$ is trivial while the non-zeroness of the $L$-function is usually proved using a hard theorem of Rohrlich.  
 
 Actually, in the process of going from the modular-symbol to the $L$-function, we perform successively three operations: first, evaluating the modular-symbol at $\{\infty\}-\{0\}$; second, restricting the obtained distribution form $\Z_p$ to $\Z_p^\ast$; third, applying the $p$-adic Mellin transform. The first two operations may not be injective.  

To come back to our question: if $f$ is $CM$ or Eisenstein, and $f_\beta$
is $\theta$-critical, 
the $L$-functions $L^\pm(f_\beta,\cdot)$ is not identically $0$ if and only if $f$ is CM, or 
$f$ is Eisenstein and $\pm=\epsilon(f_\beta)$. The last $p$-adic $L$-function 
$L^{- \epsilon(f_\beta)} (f_\beta,\cdot)$ is zero. In my sense,  this zeroness is not a defect of the definition (since the corresponding modular symbol $\Phi_{f_\beta}^{-\epsilon(f_\beta)}$ is uniquely-defined and non-zero),  but a {\it bona fide} 
property of the $L$-function, that should have some algebraic counterpart in term
of a version of the Iwasawa main conjecture.

Those results, which are easy, will be proved in the work in progress alluded above. (but see Example~\ref{exLEisenstein} below for examples of non-zero $L^{\epsilon(f_\beta)}(f_\beta,\cdot)$
when $f$ is Eisenstein)
\end{itemize}
\end{iremark}

\subsubsection{Two variables $p$-adic $L$-functions}

We note $\CC$ the eigencurve of tame level $\Gamma_1(N)$.

\begin{itheorem}\label{twovariables} Let $f_\beta$ be as in theorem~\ref{eigenspace}, and let $x$ be the point corresponding to $f_\beta$ in the eigencurve $\CC$. There exists  an affinoid 
neighborhood $V$ of $x$
in $\CC$, and for any choice of the sign $\pm$, a two-variables $p$-adic $L$-function 
$L^{\pm}(y,\sigma)$ for $y \in V$ and $\sigma$ a continuous character $\Z_p^\ast \rightarrow \C_p^\ast$, well defined up to a non-vanishing analytic function of $y$ alone. The function $L^{\pm}(y,\sigma)$ is jointly analytic in the two variables. If $y$ corresponds to a refined newform $f'_{\beta'}$, this newform is decent
and one has $L^\pm(y,\sigma)=c L^\pm(f'_{\beta'},\sigma)$ where $c$ is a non-zero scalar depending on $y$ but not on $\sigma$.
\end{itheorem}

\subsubsection{Secondary critical $p$-adic $L$-functions}

Let $f_\beta$ be as in Theorem~\ref{eigenspace}.
 Assume that $f_\beta$ is critical. The number $\dim \Symb_\Gamma^\pm(\D_k)_{(f_\beta)}$ is independent of $\pm$ and greater than one by Theorem~\ref{eigenspace}. Let us call it $e$.
We shall prove that this number is equal to the degree of ramification of the weight map $\kappa$ at the point $x$ corresponding to $f_\beta$ on the eigencurve $\CC$, and that the point $x$ is smooth.
 
 In addition to the $p$-adic $L$-function $L^{\pm}(f_\beta,\sigma)$ we can define $e-1$
 {\it secondary $L$-functions} $L^{\pm}_i(f_\beta,\sigma)$ for $i=1,\dots,e-1$, by setting,
 for all $\sigma$.
 \begin{eqnarray} \label{defsecondary} L^{\pm}_{i}(f_\beta,\sigma)= \left( \frac{\partial^i L^{\pm}(y,\sigma)}{\partial y^i} \right)_{|y=x}.
 \end{eqnarray}
 This make sense since $\CC$ is smooth of dimension $1$ at $x$. Note that extending this formula to $i=0$ gives $L_0^{\pm}(f_\beta,\sigma)=L^\pm(f_\beta,\sigma)$.
 
 Of course, the function $L_i^{\pm}(f_\beta,\sigma)$ is only defined up to a scalar and up to adding any linear combination of the $L_j^{\pm}(f_\beta,\sigma)$ for $0 \leq j \leq i-1$. In other words, this is the flag $F_0 \subset F_1 \subset \dots \subset F_{e-1}$ in the space of analytic functions of $\sigma$ that is well defined, where $F_i = \sum_{j=0}^{i} \bar \Q_p L^\pm_i(f_\beta,\sigma)$.
 
 The analytic functions $L_i^{\pm}(f_\beta,\sigma)$ for $i=0,\dots,e-1$ enjoy the following properties (which are all invariant under the indeterminacy in the definition of $L_i^{\pm}$, as the reader can easily check).

\mb
{\bf Connection with modular symbols:} There exists a unique maximal flag $F'_0 \subset F'_1 \subset \dots \subset F'_{e-1}$ (with $\dim F'_i=i+1$) in $\Symb_\Gamma(\D_k)_{(f_\beta)}$ which is mapped onto the flag $F_0 \subset \dots \subset F_{e-1}$ by the application $\phi \mapsto (\sigma \mapsto \phi(\{\infty\}-\{0\})(\sigma)$. In other words, the secondary $L$-functions of $f_\beta$ are Mellin transform of values at $\{\infty\}-\{0\}$ of generalized eigenvectors in the space of distribution-valued
modular symbols with system of eigenvalues as $f_\beta$, and those generalized eigenvectors are well-defined up to the same 
indeterminacy as the functions $L_i(f_\beta,.)$
 \mb
{\bf Order:} For all $i=0,\dots,e-1$, $L_i^{\pm}(f_\beta)$ is of order at most $v_p(\beta)=k+1$.
\mb
{\bf Interpolation:} Assume that $f$ is cuspidal. If $0 \leq i \leq e-2$, then $L_i^\pm(f_\beta,\sigma)$ satisfies the interpolation
property (\ref{intcrit}). If $i=e-1$, then $L_{e-1}^\pm(f_\beta,\sigma)$ satisfies the interpolation property 
(\ref{intnoncrit}).
\mb
{\bf Infinity of zero's} Assume that $f$ is cuspidal.
For all $i=0,\dots,e-1$, and all $\psi: \F_p^\ast \rightarrow \C_p^\ast$, the function $s \mapsto 
L_i^\pm(f_\beta,\psi \chi_s)$ has infinity many zeros.

\subsection{Method and plan of the paper}

The idea of the proof of Theorem~\ref{twovariables} is as follows: over the weight space $\WW$, or at least over an affinoid open subset $W=\sp R$ of it,
there exists a module of distributions-valued modular symbols. This module, that we note $\Symb_\Gamma(R \hotimes \D)$, or an analog to it, was already considered by other authors (e.g. \cite{greenbergstevens}): its elements give rise, essentially by Mellin transform, to two-variables functions $L(k,\sigma)$, $k$ being a variable on the open subset $W$ of the weight space $\WW$. We need to see this module, or the part of its Riesz decomposition
 corresponding to a $U_p$-slope $\leq \nu$, as a module over (the affinoid ring of a suitable affinoid subspace $\CC_{W,\nu}$ of)
the eigencurve. We have found no other way to do so than to reconstruct the eigencurve inside it, that is to use that module instead of 
Coleman's module of overconvergent modular forms in the construction of the eigencurve. It is then easy
to show using a theorem of Chenevier (\cite{chenevierJL}) that the new eigencurve thus constructed is essentially the same as the old one.
 
The second step is to study the structure of this module as a module over the affinoid $\CC_{W,\nu}$. Here we use the fundamental fact that
the eigencurve is known to be smooth at all decent classical point. Near such a point, an elementary commutative-algebra argument
shows that our module is free (of rank $2$, each of its $\pm 1$-part being of rank one). This is the result from which all the others are deduced:
one first derives that the formation of the algebra generated by the operators in $\HH$ in the modules of modular symbols commute with specialization,
which allow us to derive the structure of this algebra at a decent classical point $x$: it always has the form $\bar \Q_p[t]/(t^e)$ where $e$ is
degree of ramification of the eigencurve at $x$ over the weight $e$. This implies with a little supplementary work all result concerning 
the one-variable $p$-adic $L$-function. As for the two-variables $p$-adic $L$-function, we refer the reader to the corresponding section of the paper:~\ref{proof2}

In section 2, we prove all results on the eigencurve and the notion of refinement that we shall need later. In particular, we extend
the result of \cite{BC} that the eigencurve is smooth at Eisenstein series to decent critical cusp forms, using a simple argument due to Chenevier.

In section 3, we construct the module of modular symbols $\Symb_\Gamma(R \hotimes \D)$ (cf~\ref{famdis} and~\ref{dismodsymb}). We then reconstruct the eigencurve in it (cf~\ref{localpieces}), We also study the Mellin transform in this context (cf~\ref{mellin}),

Section 4 is devoted to the proof of the results stated in this introduction.
 
\subsection{Perspectives}

\label{perspectives}

As said earlier, the elegant and clearly exposed results of Stevens and Pollack, and the remaining open questions appearing 
clearly from their work were a strong motivation for this work. But there was a second motivation, that I wish now to explain.
As this explanation has little to do with the rest of the article, the reader may wish to skip it.
 
 From a larger point of view,
I mean for automorphic forms for larger groups, especially unitary groups, I believe it is very important, in connection to the Bloch-Kato conjecture, 
to develop a theory of $p$-adic $L$-function for automorphic forms, and of the family thereof, {\bf especially including} the critical and critical slope cases.
Let me explain why: Since 2002 a series of work by the author, Chenevier, Skinner and Urban have developed a strategy to prove 
the lower bound on the Selmer group predicted by Bloch-Kato conjecture or its $p$-adic variant. This strategy, which is part of the world of methods inspired by Ribet's proof
of the converse of Herbrand's theorem, can be summarized as follows: one picks up an automorphic form $f$ for an unitary group (or a symplectic group)
whose Galois representation is reducible (the automorphic form may be Eisenstein, or cuspidal endoscopic) and contains the trivial character, and the 
Galois representation $\rho$ of which we want to study the Selmer group. One then deform it into a family of automorphic forms whose Galois representations 
are irreducible, and thus obtain, by generalization of a famous lemma of Ribet, extensions of $1$ by $\rho$ that give the desired elements in the Selmer group.
So far, this method has only allowed to construct one extension at a time (\cite{BC},\cite{SU},\cite{BCbook}), that is to prove a lower bound of $1$
on the Selmer group\footnote{I observe that in \cite{SUicm}, the proof of a lower bound of $2$ when predicted by the sign of the Bloch-Kato conjecture
is announced. However, the proof has not been given yet (in particular, the crucial deformation argument is missing, and does not appear either in Urban's work in progress on eigenvarieties). Besides, this method does not seem to
be able to give lower bound greater than $2$ -- though a variant proposed elsewhere by the same authors does: see above.}. 
However, two different
sub-strategies have been proposed (both during the Eigensemester in 2006) to get the lower bound on the Selmer group predicted by the order of vanishing 
of the $p$-adic $L$-function of $\rho$. One, by Skinner and Urban, would use for $f$ an Eisenstein series on $\U(n+2r)$ (where $\rho$ has dimension $n$ 
and $r$ is the number of independent extensions one expect to construct, the other, by the author and Chenevier (see the final remarks in \cite{BCbook} for more
details) use for $f$ a cuspidal endoscopic non-tempered form for $U(n+2)$, and hope that the number $r$ of independent extensions will
be encoded into the singularity of the eigenvariety of $\U(n+2)$ at the point $f$, leading to the construction of $r$ independent extensions by the results of  \cite{BCbook}.
Both methods (at least it seems to me from my remembrance of Urban's explanations of his and Skinner's method at the eigensemester in 2006, since to my knowledge those explanations never appeared in print) would need
a good knowledge on the $p$-adic $L$-function on $f$, and of how it is related to the $p$-adic $L$-function of $\rho$. Such a knowledge is far from being available,
because in general little is known on higher rank $p$-adic $L$-function, but also more fundamentally because in both methods the form $f$ we work with is critical or $\theta$-critical
(and there are fundamental reasons it is so, that it would take us too far to explain here),
a case which is not well understood even for $\Gl_2$. Therefore, understanding the case of $\Gl_2$ critical $p$-adic $L$-function may be seen as a first, elementary step toward the general goal of understanding the higher rank critical $p$-adic $L$-function.

\section{Modular forms, refinements, and the eigencurve}

%\begin{itemize}
%\item[(iii)] For any quadratic extension $L/\Q$ with $L \subset \Q(\zeta_{p^3})$, $(\rhob_f)_{|G_L}$ is absolutely irreducible.
%\item[(iv)] We have $p>k$ and for the quadratic extenstion $L=\Q(\sqrt{(-1)^{(p-1)/2}p}$, the representation $(\rhob_f)_{|G_L}$ is absoultely irreducible.
%\item[(v)] The representation $\rhob_f$ is not absolutely irreducible, and its two diagonal characters are distinct after restriction to $G_{\Q(\zeta_{p^\infty})}$.
%\item[(vi)] The form $f$ is not CM, and $f$ is cuspidal or special at all prime dividing $N$ (e.g. $f$ is of trivial nebentypus and $N$ is square-free).
%\end{itemize}

\subsection{The eigencurves}

\subsubsection{The Coleman-Mazur-Buzzard eigencurve}

\label{defeigencurve}
Let $\WW$ be the usual weight space (see \S\ref{weightspace} below for more details).
We shall denote by $\CC$ the $p$-adic eigencurve of tame level $\Gamma_1(N)$ (constructed by Coleman-Mazur in the case $N=1$, $p >2$, and 
in general using similar methods by Buzzard). The space $\CC$ is equidimensional
of dimension $1$. There exists a natural  weight map $\kappa : \CC \rightarrow \WW$, and it is finite and flat. {\bf We normalize $\kappa$ in such a way that if $f_x$ is a form of weight $k+2$, then $\kappa(x)$ is $k$.}

The operators $T_l$ ($l \not | pN$), 
$U_p$ and the diamond operators, that is all elements of $\HH$ define bounded
analytic functions on $\CC$.

If $W=\sp R$ is an affinoid subspace of $\WW$, and $\nu >0$ is a real number, we shall denote by $\CC_{W,\nu}$
the open affinoid subspace of the eigencurve $\CC$ that lies above $W$ and on which one has $v_p(U_p) \leq \nu$. 
By construction, $\CC_{W,\nu}$ is the maximal spectrum of the sub-algebra of $\End_R(M^{\dag}_{W,\nu})$ generated by the image of $\HH$, where $M^\dag_{W,\nu}$
denotes Coleman's locally free of finite type $R$-modules of overconvergent modular forms of level $\Gamma=\Gamma_1(N) \cap \Gamma_0(p)$, weight in $W$ 
and slope $\leq \nu$. Let us recall that the spaces $M^{\dag}_{W,\nu}$ are defined as pieces of the Riesz decomposition for the operator $U_p$ acting on Banach's $R$-module
$M^\dag_{W}$. If $\kappa \in W(\Q_p)$ is an integral weight $k \in \N$, then there exists natural $\HH$-isomorphism
$M^\dag_W \otimes_{R,\kappa} \Q_p = M^\dag_{k+2}(\Gamma),$ where $M^\dag_{k+2}(\Gamma)$ is the space of overconvergent modular forms. This space contains as an $\HH$-submodule the $\Q_p$-space of classical modular forms $M_{k+2}(\Gamma)$.

\subsubsection{Cuspidal overconvergent modular forms and evil Eisenstein series}

\label{eisensteinseries}

The following definitions  and results are natural generalizations of  definitions and results stated in \cite{CM} in the case $N=1$.

\begin{definition} An element $f \in M^\dag_{k+2}(\Gamma)$ is called {\it overconvergent-cuspidal} if it vanishes at all cusps of $X(\Gamma)$ in the $\Gamma_1(N)$-class
of the cusp $\infty$. We call $S^\dag(k+2)(\Gamma)$ the $\HH$-submodule of 
$M^\dag_{k+2}(\Gamma)$ of overconvergent-cuspidal forms.
\end{definition}
 
 \begin{definition} An Eisenstein series in $M_{k+2}(\Gamma)$ is called 
 {\it evil} if it is overconvergent-cuspidal as an element of $M^\dag_{k+2}(\Gamma)$.
 \end{definition}
 
 Let us recall the well-known classification of Eisenstein series for the group $\Gamma_1(M)$. Let $\chi$ and $\psi$ be two primitive Dirichlet character of conductors $L$ and $R$. We assume that $\chi(-1) \psi(-1) = (-1)^k$. Let $$E_{k+2,\chi,\psi}(q)=c_0 + \sum_{m \geq 1} q^m \sum_{n|m}(\psi(n) \chi(m/n) n^{k+1})$$ where $c_0=0$ if $L >1$ and $c_0 = - B_{k,\psi}/2k$ if $L=1$. 
If $t$ is a positive integer, let $E_{k+2,\chi,\psi,t}(q)=E_{k+2,\chi,\psi}(q^t)$ excepted in the case $k=0$, $\chi=\psi=1$, where one sets $E_{2,1,1,t}=E_{2,1,1}(t)-tE_{2,1,1}(q^t)$.

\begin{prop}[Miyake, Stein] Fix an integer $M$, and a Dirichlet character $\epsilon$ of 
$(\Z/m\Z)^\ast$. 

The series $E_{k+2,\chi,\psi,t}(q)$ are modular forms of level $\Gamma_1(M)$ and character $\epsilon$ for all positives integers 
$L$, $R$, $t$ such that $LRt  | M$ and all primitive Dirichlet character $\chi$ of conductor $L$ and $\psi$ of conductor $M$, satisfying $\chi \psi = \epsilon$
(and $t>1$ in the case $k=0$, $\chi=\psi=1$) and moreover they form a basis of the space of Eisenstein series in $M_{k+2}(\Gamma_1(M),\epsilon)$. 

For every prime $l$ not dividing $M$, we have 
$$T_l E_{k+2,\chi,\psi,t}= (\chi(l) + \psi(l) l^{k+1}) E_{k+2,\chi,\psi,t}.$$
\end{prop}
\begin{pf} See Miyake (\cite{miyake}) for the computations leading to those results, Stein (\cite{Stein}) for the results stated as here. (The eigenvalues for $T_l$ are given in \cite{Stein} only in the case $t=1$, but the general case follows immediately using the relations between Hecke eigenvalues and Fourier coefficients for non-necessarily normalized modular forms)
\end{pf}
From this description it follows easily that 
\begin{cor} The Eisenstein series in $M_{k+2}(\Gamma_1(N))$ that are new in the sense of the introduction
are exactly the following: the Eisenstein series $E_{k+2,\psi,\chi,1}$ with $M=LR$ (excepted of course $E_{2,1,1,1}$ which is not even a modular form). The Eisenstein series $E_{2,1,1,M}$'s when $M$ is prime. Those series generate their eigenspaces for the Hecke operators $T_l$, $l$ does not divide $M$, 
in $M_{k+2}(\Gamma_1(M),\psi \chi)$. 
\end{cor}

We apply those results to the situation of our article:
\begin{cor} \label{coreis} Let $f$ be an Eisenstein series of weight $k+2$ and level 
$\Gamma_1(N)$, eigenform for the operator $T_p$, and let $f_\alpha$ and $f_\beta$ be its two refinement (forms for $\Gamma=\Gamma_1(N) \cap \Gamma_0(p)$. 
\begin{itemize}
\item[(i)] The valuation of $\alpha$ and $\beta$ are, in some order $0$ and $k+1$.
\item[(ii)] If $f$ is new, the $\HH$-eigenspace of $f_\alpha$ and $f_\beta$ are of dimension $1$.
\item[(iii)] The form $f_\beta$ is evil if and only if it is of critical slope, that is if 
$v_p(\beta)=k+1$. 
\end{itemize}
\end{cor}
\begin{pf} By the proposition, $f$ lies in the same $T_p$ eigenspace as
some Eisenstein series $E_{k+2,\chi,\psi,t}$, so if $T_p f = a_p f$ we have
$a_p = \chi(p) + \psi(p) p^{k-1}$. Therefore, the roots $\alpha$ and $\beta$ of the polynomial (\ref{hp}) are $\chi(p)$ and $\psi(p) p^{k+1}$, and (i) follows.

Let us prove (ii). If $f$ is new then $f=E_{k+2,\psi,\chi,1}(q)$ with $\psi$ and $\chi$ primitive of conductor $L$ and $R$,  and $N=LR$ by the corollary (we leave the exceptional case $E_{2,1,1,M}$ to the reader).
Applying the Proposition to $M=Np$ show that the space of Eisenstein series for 
$\Gamma_1(N) \cap \Gamma_0(p)$ and character $\chi \psi$ (that is
the space of Eisenstein series for  $\Gamma_1(Np)$ and character $\chi \psi$ seen as a character of conductor $Np$) has dimension $2$,
and is generated by $E_{k+2,\chi,\psi,1}(q)=f(z)$ and 
$E_{k+2,\chi,\psi,p}(z)=f(pz)$. Since $U_p$ has two distinct
eingevalues $\alpha$ and $\beta$ on that space whose eigenvectors are respectively
$f(z)-\beta f(pz)$ and $f(z)-\alpha f(pz)$, the result follows.

To prove (iii), we note that the number of cusps for $\Gamma$ that are
 $\Gamma_1(N)$-equivalent to $\infty$ is half the number of cusps for $\Gamma=\Gamma_1(N) \cap \Gamma_0(p)$, so the dimension of the space of evil Eisenstein series is half the dimension of the space of all
Eisenstein series (putting the exceptional case aside, where the reader can take care of it), as is obviously the space generated by critical slope Eisenstein series $f_\beta$. Therefore, it is enough
to show that critical slope Eisenstein series are evil, as the converse would follow from equality of dimensions. So let $f$ be an Eisenstein for $\Gamma_1(N)$, and we can
assume that $f=E_{k+2,\psi,\chi}$ (leaving once again the exceptional case to the reader). If $v_p(\beta)=k+1$, then $f_\beta(z)=
f(z)-\alpha f(p z)$. If $\chi=1$, i.e. if $L=1$, then $\alpha=\chi(p)=1$
so $f_\beta(\infty)=f(\infty)-f(\infty)=0$. If $\chi \neq 1$, i.e. $L \neq 1$, we have $f(\infty)=c_0=0$ by the Proposition, so $f_\beta(\infty)=0$. We thus see that in both cases, 
$f_\beta$ vanishes at the cusp $\infty$. Applying the same reasoning to $f_{|\gamma}$ for $\gamma \in \Gamma_1(N)$ shows that $f_\beta$ vanishes at all cusps
that are $\Gamma_1(N)$-equivalent to $\infty$. Therefore, $f_\beta$ is evil.
\end{pf}

Let us also note the following
\begin{cor}[Control theorem for cuspidal forms]   \label{controlcusp}
 $$S_{k+2}^\dag(\Gamma)^{< k+1} = S_{k+2}(\Gamma)^{<k+1}.$$
\end{cor}
\begin{pf} By Coleman's control theorem \cite{coleman}, 
any form $f \in S_{k+2}^\dag(\Gamma)^{<k+1}$ is classical. We can 
assume that $f$ is an eigenform, so the question
is whether it is cuspidal or Eisenstein. But if $f$ was Eisenstein, it would be evil
by definition, so of critical slope $k+1$ by point (iii) of Corollary~\ref{coreis}, which is absurd since it is of slope $<k+1$.
\end{pf}

 \subsubsection{The cuspidal eigencurve}

\label{cuspidaleigencurve}

Mimicking Coleman's construction, 
It is easy to construct a cuspidal analog $S^\dag_W$ of $M^\dag_W$, which is a Banach $R$-module (even satisfying property (Pr) of \cite{buzzard})
 with an action of $\HH$ (and a $\HH$-equivariant embeding $S^\dag_W \hookrightarrow M^\dag_W$) such that $U_p$ acts compactly, and with natural $\HH$-isomorphisms $S^\dag_W \otimes_{R,\kappa} \Q_p = S^\dag_{k+2}(\Gamma),$ if $\kappa \in W(\Q_p)$ is an integral weight $k$. The general eigenvariety machine (\cite{buzzard}) then produces a {\it cuspidal eigencurve} $\CC^0$ of tame level 
 $\Gamma_1(N)$, and by construction $\CC^0$ is equidimensional of dimension $1$, provided with a map $\kappa : \CC^0 \rightarrow \WW$ which is locally finite and flat,
and is naturally identified with a closed subset of $\CC$.

\subsubsection{Degree of the weight map and overconvergent modular forms}

Let us prove for further reference two results about classical and overconvergent modular newforms:
\begin{lemma} \label{newnew} 
Let $x$ be a point of $\CC$ corresponding to a refined newform. Then in a neighborhood
of $\CC$ in $x$ all classical points also correspond to refined {\bf newforms}.
\end{lemma}
\begin{pf} For $N'$ a proper divisor of $N$ let $\CC_{N'}$ denote the eigencurve of tame level $\Gamma_1(N')$ constructed with the same algebra $\HH$ we have considered so far, that is with no operators $T_l$, for $l$ a divisor of $N/N'$. By construction, there is closed immersion $\CC_{N'} \hookrightarrow \CC$. Since the image of no $\CC_{N'}$ contains $x$ (since $x$ is new), and since
there is only a finite number of proper divisors $N'$ of $N$,
the same is true in a neighborhood of $x$, and this implies the lemma.
\end{pf}

\begin{remark} In \cite{newton} a stronger result is proved, namely that
the Zariski-closure of {\bf new} classical points in a suitable neighborhood of $x$
is of equidimension $1$.  We shall not need this result, and actually it will be a
consequence of this lemma for decent classical $x$'s once we have proved 
Theorem~\ref{smooth}. But Newton's result is much more general since it applies as well
to non-classical points, in particular to points where two components of the eigencurve 
meet (and of course, the eigencurve is non smooth at such points).
\end{remark}
 
\begin{lemma} \label{lemmanewform} Let $f$ be a newform of level $\Gamma_1(N)$ and weight $k+2$, and let $\alpha$ and $\beta$ be the two roots
of the quadratic polynomial~(\ref{hp}). Let $f_\beta$ be one of its refinements, and $x$ the corresponding point in $\CC$, that we identify with the system of eigenvalues $x: \HH \rightarrow \bar \Q_p$ attached to $f_\beta$.
\begin{itemize}
\item[(i)] The generalized eigenspace $M_{k+2}(\Gamma)_{(x)}$ has dimension $1$.
\item[(ii)] The generalized eigenspace $M_{k+2}^\dag(\Gamma)_{(x)}$ has dimension $e$, where $e$ is the length of the affine ring of the scheme-theoretic 
fiber of $\kappa:\CC \rightarrow \WW$ at $x$.
\end{itemize}
If we assume in addition that $f$ is cuspidal, or that $f$ is Eisenstein and $v_p(\beta)>0$, then: 
\begin{itemize}
\item[(iii)] The point $x$ belongs to $\CC^0$ and the generalized eigenspace $S_{k+2}^\dag(\Gamma)_{(x)}$ has dimension $e^0$, where $e^0$ 
is the length of the affine ring of the scheme-theoretic fiber of $\kappa:\CC^0 \rightarrow \WW$ at $x$.  (We shall see in Cor.~\ref{corsmooth} that very often $e=e_0$.)
\end{itemize}
\end{lemma}
\begin{pf}
The first result is a well-known consequence of the theory of newforms in the cuspidal case, and of Corollary~\ref{coreis}(ii) in the Eisenstein case.

Let us prove the second result. First, we note that it reduces to (i) in the easy case where $\kappa$ is \'etale at $x$ and $v_p(U_p(x))<k+1$, for in this case, 
$e=1$, and $M_{k+2}^\dag(\Gamma)_{(x)} = M_{k+2}(\Gamma)_{(x)}$ by Coleman's control theorem. 

For a general $x$ we can assume, shrinking $W$ if necessary, that 
 $x$ is the unique geometric point above $k$ in the connected component $X$ of $x$ in $\CC_{W,\nu}$, and that $v_p(U_p)$ is constant on $X$. Let us denote by $M_X^\dag$ the sub-module of $M_W^\dag$ corresponding to the component $X$. This is a direct summand of $M_W^\dag$, hence it is a finite projective module over $R$. One has
by definition for every $k' \in \N \cap W(\Q_p)$ a natural isomorphism of $\HH$-module 
\begin{eqnarray} 
\label{MX} M_X^\dag \otimes_{R,k'} \bar \Q_p = \oplus_{y \in X, \kappa(y)=k'} M^\dag_{k'+2}(\Gamma)_{(y)}.
\end{eqnarray} 

The degree of the flat map  $\kappa_{|X}: X \rightarrow W$ is $e$ by definition of $e$, so we have, for all points $k' \in \N \cap W(\Q_p)$ where $\kappa_{|X}$ is \'etale, $e=|\kappa_{|X}^{-1}(k')(\bar \Q_p)|$. If moreover
$k'+1 > v_p(U_p(x))$, then for any point $y$ in $\kappa_{|X}^{-1}(k')(\bar \Q_p)$, 
such that $y$ corresponds to a refined newform, one has $\dim M_{k'+2}^\dag(\Gamma)_{(y)}=1$ by the easy case
treated above. Therefore, by (\ref{MX}), for such a $k'$,
 the space $M^\dag_X \otimes_{R,k'} \bar \Q_p$ has dimension $e$. Since, in view of Lemma~\ref{newnew}, 
we can find  integers $k'$'s arbitrary closed in $W$ to $k$ satisfying the above conditions, we deduce that the rank of $M_X^\dag$ is $e$. 
By applying $(\ref{MX})$ at $k$, we see that  $\dim M_{k+2}^\dag(\Gamma)_{(x)}=e$ since $x$ is the unique geometric point of $X$ above $k$.

As for (iii), we note that it is clear that $x$ belongs to $\CC^0$ if $f$ is cuspidal; if $f$ is Eisenstein, the hypothesis $v_p(\beta) > 0$ implies that $f_\beta \in S^\dag_{k+2}(\Gamma)$ by Corollary~\ref{coreis}(iii). Therefore $x$ belongs to $\CC^0$. 
The rest of the proof of (iii) is exactly as in (ii) replacing $\CC$ by $\CC^0$ and $M^\dag$ by $S^\dag$ everywhere. 
\end{pf}
\begin{remark} \label{remarknewform} The proof of this lemma is close to the proof of \cite[Prop. 1(b)]{BCsmooth}.
Note that while \cite{BCsmooth} was written with the assumption that $N=1$ (the only case the eigencurve was constructed 
at that time), the arguments and results of this paper are still correct for any $N$ without any change if we simply use 
Lemma~\ref{lemmanewform} instead of~\cite[Prop. 1(b)]{BCsmooth} and work with Eisenstein series that are newforms.
In what follows, we shall use without further comment the results of \cite{BCsmooth} for any level $N$.
\end{remark}

\subsection{Refinements}

\subsubsection{The preferred Galois representation}

As in the introduction, let $f=\sum a_n q^n$ be a modular newform (cuspidal or eisenstein) of weight $k+2$ and level $\Gamma_1(N)$,
  \label{preferred}
\begin{lemma}
There exists a unique (up to isomorphism) Galois representation $\rho_f : G_\Q \rightarrow \Gl_2(\bar \Q_p)$ such that
\begin{itemize} 
\item[(i)] The representation $\rho_f$ is unramified outside $Np$ and 
satisfies the Eichler-Shimura relation $\tr(\rho_f(\Frob_l)) = a_l$ for all $l$ prime to $Np$.
\item[(ii)] The restriction $(\rho_f)_{|D_p}$ is crystalline at $p$.
\item[(iii)] The representation $\rho_f$ is indecomposable.
\end{itemize}
 Moreover, the Hodge-Tate weights of  $(\rho_f)_{|D_p}$ are $0$ and $k+1$ and the eigenvalues of the crystalline Frobenius $\varphi$ 
on $D_\crys((\rho_f)_{|D_p})$ are the roots $\alpha$ and $\beta$ of the quadratic polynomial (\ref{hp}). 
\end{lemma}
\begin{pf} 
In the case where $f$ is cuspidal, the existence and uniqueness of $\rho_f$ satisfying (1) are well-known, 
as well as its properties at $p$, and $\rho_f$ is irreducible in addition. In the case where $f$ is Eisenstein, 
what is well-known is the existence of a decomposable representation $\chi_1 \oplus \chi_2$ satisfying (i), crystalline at $p$, with say 
$\chi_1$ of Hodge-Tate weight $0$ and $\chi_2$ of Hodge-Tate weight $k+1$. In this case the quotient $\chi_1 \chi_2^{-1}$ is $\omega^{k+1} \epsilon$
(where $\epsilon$ is the Dirichlet character of $(\Z/N\Z)^\ast$ defined by $\langle a \rangle f = \epsilon(a) f$) and by the classification of Eisenstein series, we have $\epsilon(-1)=1$, and moreover $\epsilon \neq 1$ when $k=0$ 
(since ``$E_2$ does not exist''). The conjecture of Bloch-Kato, known by the work of Soul\'e for characters of $G_\Q$, says that
$\dim \Ext^1_f(\chi_2,\chi_1)=1$ and $\dim \Ext^1_f(\chi_1,\chi_2)=0$, where $\Ext^1_f$ parametrizes the extensions that are crystalline at $p$.
Therefore, there exists a unique representation $\rho_f$ satisfying (i), (ii) and (iii) up to isomorphism, namely the extension of $\chi_2$ by $\chi_1$ defined by any non-zero element in $\Ext^1_f(\chi_2,\chi_2)$.
\end{pf}

We call $\rho_f$ the {\it preferred representation} attached to $f$. Note that since $\alpha \neq \beta$ (by assumption)
the choice of a refinement $f_\alpha$ or $f_\beta$ of $f$ is equivalent to the choice of a $\varphi$-stable line in 
$D_\crys((\rho_f)_{|D_p})$, namely the eigenspace of the crystalline Frobenius 
 $\varphi$ of eigenvalue $\alpha$ or $\beta$ respectively. This line is called a {\it refinement} of $\rho_f$, and is the refinement of $\rho_f$ attached to the given refinement of $f$.

\subsubsection{Refinements and the preferred Galois representation}

\label{refinements}
\label{refinement}

Note that we have obviously $v_p(\alpha) + v_p(\beta) = k+1$ (recall that our forms are of weight $k+2$, not $k$), and $0 \leq v_p(\alpha), v_p(\beta) \leq k+1$. We shall say that a refined modular form $f_\beta$ is
{\it of critical slope} if $v_p(\beta)=k+1$. Then of course, for a given $f$ as above, there is at least one
of the two refinements $f_\beta$ and $f_\alpha$ of $f$ which is of non-critical slope, and exactly one
if and only if $v_p(a_p) = 0$ (in which case $f$ is called {\it ordinary}).
  
\begin{prop}\label{propcrit} For a refined form $f_\beta$, the following are equivalent.
\begin{itemize} 

\item[(i)] There exist {\bf non-classical} overconvergent modular forms of weight $k+2$ which are {\bf generalized} eigenvectors of $\HH$ with the same eigenvalues as $f_\beta$. In other words, the inclusion $M_{k+2}(\Gamma)_{(f_\beta)} \subset M_{k+2}^\dag(\Gamma)_{(f_\beta)}$ is strict.
\item[(i')] We have $\dim M_{k+2}^\dag(\Gamma)_{(f_\beta)} \geq 2$.
\item[(ii)] We have $v_p(\beta)>0$ and the eigencurve $\CC$ is {\bf not \'etale} over the weight space at the point $x$ corresponding to $f_\beta$ (see below)
\item[(iii)] We have $v_p(\beta)>0$ and the representation $(\rho_f)_{|D_p}$ is the direct sum of two characters.
\item[(iv)] The refinement of $\rho_f$ attached to $f_\beta$ is critical in the sense of \cite{BCbook}, that is the line of eigenvalue $\beta$
in $D_\crys((\rho_f)_{|D_p})$ is the same as the line defined by the weight filtration.
\end{itemize}
Moreover if those property holds, $v_p(\beta)=k+1$, that is $f_\beta$ is of critical slope.

In the case $f$ is cuspidal, those properties are also equivalent to 
\begin{itemize}
\item[(v)] The modular form $f_\beta$ is in the image by the operator $\theta_k: M_{-k}^\dag(\Gamma) \rightarrow  M_{k+2}^\dag(\Gamma)$ which acts on $q$-expansion by $(q\frac{d}{dq})^{k+1}$.
\end{itemize}
In the case $f$ is Eisenstein, {\rm{(v)}} is always true.
\end{prop}
\begin{pf}

All this is well known and we just give references or short indications of proof.

The equivalences between (i), (i')  and (ii) follow directly from Lemma~\ref{lemmanewform}.

The equivalence between (iii) and (iv) follows easily from the weak admissibility of $D_\crys((\rho_f)_{|D_p}$, as does the assertion after the moreover:  see \cite[Remark]{BCbook} for details.

To finish the proof we separate the Eisenstein and cuspidal cases.

In the Eisenstein case, we only have to show that (i) is equivalent to (iii) to complete the proof of the equivalence of (i) to (iv). But this is a theorem of  the author
and Chenevier (\cite[Th\'eor\`eme 3]{BCsmooth}). We also have to check that (v) always hold.
But the limit case  of Coleman's control theorem (cf. \cite[Corollary 7.2.2]{coleman}) states that the cokernel 
of $\theta_k$ in $(M_{k+2}(\Gamma)^\dag)_{(x)}$ is isomorphic to $S_{k+2}(\Gamma)_{(x)}$, that is $0$,
or in other words, $\theta_k$ is surjective onto $(M_{k+2}(\Gamma)^\dag)_{(x)}$, so obviously $f_\beta$ is in its image. 

In the cuspidal case, we note that the equivalence between (v) and (iii) is a result of Breuil and Emerton (\cite[Theorem 2]{BE}), so it only remains to prove that (v) is equivalent to (i') to complete the proof. 
But again by the limit case of Coleman's control theorem, the image of $\theta_k$ has codimension $1$
in $M_{k+2}(\Gamma)^{\dag}_{(x)}$, so (i') is equivalent to the assertion that the image of $\theta_k$ is not $0$. Therefore it is clear that (v) implies (i'). Conversely, since $\theta_k$ is Hecke-covariant, if it is non-zero there is a non-zero eigenvector in its image. By looking at $q$-expansions (as in \cite{stevenspollack2}, proof of Theorem 7.2), this eigenvector has to be $f_\beta$, which proves (v).
\end{pf}

\begin{definition} We say that $f_\beta$ is {\it critical} if conditions (i) to (iv) are satisfied and {\it $\theta$-critical} if condition (v) is
satisfied.
\end{definition}

\subsubsection{Classification of critical and critical-slope refined modular form}

\begin{prop} \label{criticalclass} Let $f$ be a newform as above. Then there exists a refinement $f_\beta$ of $f$ of critical slope in and only in the following cases (and in those cases there exist exactly one critical-slope refinement)
\begin{itemize}
\item[(a)] the newform $f$ is cuspidal, not CM, and ordinary at $p$ (that is $v_p(a_p)=0$).
\item[(b)] The newform $f$ is cuspidal CM, for a quadratic imaginary field $K$ where $p$ is split.
\item[(c)] The newform is Eisenstein.
\end{itemize}
In any of those case, let $f_\beta$ be the critical-slope refinement. Then
\begin{itemize}
\item[-] In case (a), it is expected, but not known in general, that $f_\beta$ is non-critical (hence not $\theta$-critical).
\item[-] In case (b), $f_\beta$ is always critical (hence $\theta$-critical).
\item[-] In case (c), $f_\beta$ is always $\theta$-critical, and it is expected, but not known in general, that $f_\beta$ is non-critical. It is known however that $f_\beta$ is non-critical when $p$ is a regular prime, or for any given $p$ and $N$ that at most a finite number of $f_\beta$ are critical.
\end{itemize}
\end{prop}
\begin{pf}
Since $v_p(\beta)+v_p(\alpha)=k+1$, and $v_p(\beta) \geq 0$, $v_p(\alpha) \geq 0$, if one of those valuation is $k+1$, then the other is $0$,
and therefore we have $v_p(a_p)v_p(-\beta - \alpha)=0$. This means that $f$ is ordinary at $p$. Of course, this is automatic for an Eisenstein series,
and equivalent to $p$ being split in $K$ for a CM form of quadratic imaginary field $K$. This proves the first set of assertions.

As for the second, it is clear that if $f$ is cuspidal CM, $(\rho_f)_{|D_p}$ is the direct sum of two characters, so $f_\beta$ is critical by 
Prop.~\ref{propcrit}(iii). If $f$ is cuspidal not CM, a well-known conjecture asserts that $(\rho_f)_{|D_p}$ is indecomposable,
and if this were true, the non-criticality of $f_\beta$ would follow. As for Eisenstein series, see~\cite[Th\'eor\`eme 3]{BCbook} and the following remark.
\end{pf}

\subsubsection{Decent refined modular forms}

\label{decent}

We recall

\begin{prop} Let $f$ be a cuspidal newform of level $\Gamma_1(N)$ and weight $k+2$ and $p$ a prime number not dividing $N$. Then if any of the following holds
\begin{itemize}
\item[(i)] For any quadratic extension $L/\Q$ with $L \subset \Q(\zeta_{p^3})$, $(\rhob_f)_{|G_L}$ is absolutely irreducible.
\item[(ii)] The representation $\rhob_f$ is not absolutely irreducible, and its two diagonal characters are distinct after restriction to $G_{\Q(\zeta_{p^\infty})}$.
\item[(iii)] The form $f$ is not CM, and $f$ is cuspidal or special at all primes dividing $N$ (e.g. $f$ is of trivial nebentypus and $N$ is square-free).
\item[(iv)] The form $f$ is CM.
\end{itemize}
we have $H^1_f(\G_\Q,\ad \rho_f)=0$, and in particular $f_\beta$ is decent for any choice of the root $\beta$ of (\ref{hp}).
\end{prop}
\begin{pf} cf. \cite{kisingeom} in cases (i) and (ii), and \cite{weston} in case (iii). In case (iv), it suffices to show by inflation-restriction that $H^1_f(G_K,\ad \rho_f)=0$ where $K$ is quadratic imaginary field attached to $f$. But $(\rho_f)_{G_K}$ is the sum of two characters $\chi$ and $\chi^c$(after extending the scalar to $\bar \Q_p$), conjugate under the outer automorphism of $G_K$ defined by the non-trivial element in $\Gal(K/\Q)$. Therefore,
$H^1_f(G_K,\ad \rho_f)=H^1_f(G_K,\chi \otimes \chi^{-1}) \oplus H^1_f(G_K.\chi \otimes (\chi^{c})^{-1}) \oplus H^1_f(G_K,\chi^c \otimes \chi^{-1})
\oplus H^1_f(G_K,\chi^c \otimes (\chi^c)^{-1})$. The first and last terms are $H^1_f(G_K,\bar \Q_p)$ which are $0$ by Class Field Theory.
The second and third terms are $H^1_f(G_K,\psi)$ where $\psi$ is a character of motivic weight $0$, and their nullity result by standard arguments from
the main conjecture for imaginary quadratic fields proved by Rubin (\cite{R}).  \end{pf}

\subsection{Smoothness of the eigencurve and consequences}

The proof  of  the following theorem is due to G. Chenevier (personal communication).
In the case of an Eisenstein series, it was proved earlier in \cite{BCsmooth}
with a partially different method (which gives more precise information). We call $G_\Q$ the Galois group of $\bar \Q$ over $\Q$ and $D_p=G_{\Q_p}$ its decomposition subgroup at $p$.

\begin{theorem} \label{smooth} Let $f$ be a decent newform of weight $k+2$ and level $\Gamma_1(N)$ and $f_\beta$ be a refinement of $f$. If $f$ is Eisenstein, assume that $v_p(\beta) \neq 0$.
let $x$ be the point on $\CC^0$ (see Lemma~\ref{lemmanewform}) corresponding to $f$. 
Then $\CC$ and $\CC^0$ are smooth at $x$.
\end{theorem} 
\begin{pf}
It suffices to prove the result for $\CC$, since $x \in \CC^0$ and $\CC^0$ is equidimensional of dimension $1$. 

When $f$ is Eisenstein, this result is the the first main theorem of \cite{BCsmooth}.
If $f_\beta$ is non-critical, then the eigencurve is \'etale at $x$ by Lemma~\ref{lemmanewform}, so it is smooth.
Therefore, since $f$ is decent, we can assume that $f$ is cuspidal, and that $f_\beta$ is critical. Thus 
it satisfies  $$H^1_f(\G_\Q,\ad \rho_f)=0.$$

Since $f$ is cuspidal, the Galois representation $\rho_f$ of $G_\Q$ attached to it is irreducible. Therefore, by a well-known theorem of Rouquier and Nyssen, 
the Galois pseudocharacter $T:G_\Q \rightarrow \TT_x$, where $\TT_x$ is the completed local ring of the eigencurve $\CC$ at the point $x$ corresponding 
to $f_\beta$ (after scalar extension  to $\bar \Q_p$)
is the trace of a unique representation $\rho: G_\Q \rightarrow \Gl_2(\TT_x)$ whose residual representation is $\rho_f$. We have
normalized $T$ and $\rho$ so that $0$ is a constant Hodge-Tate-Sen weight. In this case, Kisin's theorem (see \cite{kisin} or \cite[chapter 3]{BCbook}) states that for all
ideal $I$ of cofinite length of $\TT_x$, we have $D_\crys((\rho \otimes \TT_x/I)_{|D_p})^{\varphi=\tilde \beta}$ is free of rank one for
some $\tilde \beta \in \TT_x / I$ lifting $\beta$. 

We consider the following deformation problem: 
for all Artinian local algebra $A$ with residue field $\bar \Q_p$, we define $D(A)$ as the set
of strict isomorphism classes of representations $\rho_A: G_\Q \rightarrow \Gl_2(A)$ such that $\rho_A$ has a constant weight equal to $0$
and $D_\crys((\rho_A)_{|D_P})^{\varphi=\tilde \beta}$ is free of rank one over $A$ for some $\tilde \beta$ lifting $\beta$.
By~\cite[Prop. 8.7]{kisin} (using the fact that $\beta \neq \alpha$ and that $0 \neq k+1$), 
we see that $D$ is pro-representable, say by a
complete Noetherian local ring $R$. Since $\rho \otimes \TT_x/I$ is for all cofinite length ideal $I$ in $\TT_x$ an element of $D(\TT_x/I)$, $\rho$ define 
a morphism of algebras $R \rightarrow \TT_x$. A standard argument shows that this morphism is surjective.
Since $\TT_x$ has Krull dimension $1$, if we prove that the tangent space of $R$ has dimension at most $1$, it would follow that the map $R \rightarrow \TT_x$ is an isomorphism and that $R$ is a regular ring of dimension $1$. This would complete our proof.

The tangent space of $R$, $t_d := D(\bar \Q_p[\vareps])$, lies inside the tangent space of the deformation ring of $\rho_f$ without local condition, which is 
canonically identified with $H^1(G_\Q,\ad \rho_f)$. Since $H^1_f(G_\Q,\ad \rho_f)=0$, this space injects into $H^1_{\!/\!f}(D_p,\ad \rho_f)$ (where $H^1_{\!/\!f}$ means $H^1/H^1_f$) and so does $t_v$. Since $f_\beta$ is critical, $(\rho_f)_{|D_p}$ is by Prop.~\ref{propcrit} the direct sum of two characters, $\chi_1$ and $\chi_2$. Both are crystalline, and say $\chi_1$ has weight $0$ while $\chi_2$ has weight $k+1$. Since $v_p(\beta)=k+1$, $\beta$ is the eigenvalue of the crystalline Frobenius on $D_\crys(\chi_2)$. 
Then we compute: 
$$H^1_{\!/\!f}(D_p,\ad \rho_f) = H^1_{\!/\!f}(D_p,\chi_1 \chi_1^{-1}) \oplus H^1_{\!/\!f}(D_p,\chi_2\chi_1^{-1}) \oplus H^1_{\!/\!f}(D_p,\chi_1 \chi_2^{-1}) \oplus 
H^1_{\!/\!f}(D_p,\chi_2 \chi_2^{-1}).$$
The condition on the constant weight $0$ in our deformation problem $D$ implies that the image of $t_D$ in the first factor is $0$. The condition on 
$D_\crys(-)^{\varphi=\tilde \beta}$ implies that the image of $t_d$ in the third and fourth factors are $0$ (the third factor is $0$ anyway). Therefore
$t_D$ injects in $H^1_{\!/\!f}(D_p,\chi_2 \chi_1^{-1})=H^1(D_p,\chi_2\chi_1^{-1})$. Since $\chi_2\chi_1^{-1}$ is not trivial (look at the weight) nor the cyclotomic character (this would implies that $f$ is of weight $2$ and level $1$, and there are no such modular form), local Tate duality and Euler characteristic formula implies
that $H^1(D_p,\chi_2\chi_1^{-1})$ has dimension $1$. Therefore, $t_d$ has dimension at most one, which is what remained to prove.
\end{pf}
  
\begin{cor} \label{corsmooth} With the notation and hypotheses of the theorem, $\CC$ and $\CC^0$ are locally isomorphic at $x$ and one has 
$ S^\dag_{k+2}(\Gamma)_{(f_\beta)} =  M^\dag_{k+2}(\Gamma)_{(f_\beta)}$
\end{cor}
\begin{pf} The first assertion is clear and the second follows from Lemma~\ref{lemmanewform}(ii) and (iii) since with the notations
of this lemma we clearly have $e=e^0$.
\end{pf} 

\section{Construction of the eigencurve through distributions-valued modular symbols}

In this section we use ideas of Stevens to give a construction of two 
versions $\CC^+$ and $\CC^-$ of the eigencurve through overconvergent modular symbols. 
We then use a result of Chenevier to show that these eigencurves are actually 
canonically embedded in the Coleman-Mazur eigencurve $\CC$ and contain  its cuspidal locus $\CC^0$.

To be more precise we shall only construct the standard local pieces 
$\CC^\pm_{W,\nu}$  of the eigencurves, parametrizing modular symbols of weight belonging to a fixed affinoid subset $W$ 
of the weight space $\WW$, and of slope bounded above by $\nu$. This is these local pieces that we shall prove, using Chenevier's result, to be canonically isomorphic
to the corresponding local pieces $\CC_{W,\nu}$ and $\CC^0_{W,\nu}$ of the usual eigencurve.

\subsection{Families of distributions} 

\label{famdis}

\subsubsection{The weight space}
\label{weightspace}
Let $\WW = \Hom_{\cont,\text{group}}(\Z_p^\ast,\G_m)$ be the standard one-dimensional {\it weight space}. As is well known, $\WW$ has a structure of rigid analytic space over $\Q_p$, and is the union of $p-1$ open balls $\WW_i$, 
each of then corresponding to characters whose restriction to the cyclic subgroup $\mu_{p-1}$ of $(p-1)$-th roots of unity in $\Z_p^\ast$ is the raising to the power $i$, for $i=1,2,\dots,p-1$. If $L$ is any normed extension of $\Q_p$, a point $w \in \WW(L)$ corresponds by definition to a continuous morphism of groups 
$\tw :  \Z_p^\ast \rightarrow L^\ast$.

If $k \in \N$, then the morphism $\tilde k (t)=t^k$ for all $t \in \Z_p^\ast$ is a character $\Z_p^\ast \rightarrow \Q_p^\ast$, 
and we denote by $k$ the corresponding point in $\WW$. Thus $\N \subset \WW(\Q_p)$, and we call points in $\N$ of $\WW(\Q_p)$ {\it integral weights}.
 
We fix once and for all an open affinoid subspace $W =\sp R$ of some $\WW_i$.
The natural immersion $\sp R = W \hookrightarrow \WW$ defines a canonical
element $K \in \WW(R)$, that is a continuous morphism of groups
$$K : \Z_p^\ast \rightarrow  R^\ast,$$
whose restriction to the subgroup of $(p-1)$-th roots of unity is the raising to the power $i$.
If $w \in W(L)$, where $L$ is any normed $\Q_p$-extension, 
we can see $w$ as a ring morphism $R \rightarrow L$, and we have formally
\begin{eqnarray} \label{Kw} \tw := w \circ K.\end{eqnarray}
 
\subsubsection{Spaces of distributions}

For $r \in |\Q_p^\ast| = p^\Z$, we set, following Stevens,
$$B[\Z_p,r] = \{z \in \C_p, \exists a \in \Z_p, |z-a| \leq r\}.$$
This has a structure of an affinoid space over $\Q_p$, and we denote
by $\A[r]$ the ring of affinoid functions on $B[\Z_p,r]$. For $f \in \A[r]$, we set
$\NN{f}_r=\sup_{z \in B[\Z_p,r]} |f(z)|$. This norm makes $\A[r]$ a $\Q_p$-Banach algebra.

We denote by $\D[r]$ the continuous dual $\Hom(\A[r],\Q_p)$ of $\A[r]$. 
It is also a Banach space over $\Q_p$ for the norm $\NN{\mu}_r = \sup_{f \in \A[r]} \frac{|\mu(f)|}{\NN{f}_r}.$

There are natural restriction maps $\A[r_1] \rightarrow \A[r_2]$ for $r_1 > r_2$, that are well-known to be compact. So are their 
transposed maps $\D[r_2] \rightarrow \D[r_1]$.

\subsubsection{Modules of distributions}

We now consider the Banach $R$-modules $R \hotimes_{\Q_p} \A[r]$ and $R \hotimes_{\Q_p} \D[r]$. 
\begin{lemma} \label{ON}
\begin{itemize}
\item[(i)] Those modules are orthonormalizable Banach $R$-modules.
\item[(ii)] For $r_1>r_2$, the natural maps $R \hotimes \A[r_1] \rightarrow R \hotimes \A[r_2]$ are compact.
and $R \hotimes \D[r_2] \rightarrow R \hotimes \D[r_1]$  
\item[(iii)] We have $\Hom_R (R \hotimes_{\Q_p} \A[r] , R) = R \hotimes_{\Q_p} \D[r]$.
\end{itemize} 
\end{lemma}
\begin{pf} 
Since the norm of every element of $\A[r]$ is an element of $|\Q_p^\ast|$ (by our hypothesis that $r \in |\Q_p^\ast|$), and since the same result
holds, as a consequence, for $\D[r]$, those Banach spaces are orthonormalizable
by \cite[Prop. 1]{serre}. Therefore, $R \hotimes \A[r]$ and $R \hotimes \D[r]$
are  orthonormalizable Banach $R$-modules by \cite[Lemma 2.8]{buzzard}, which proves (i). The point (ii) follows form \cite[Corollary 2.9]{buzzard}.
Point (iii) is left to the reader.\end{pf} 

\subsubsection{Actions of $\Sigma_0(p)$}

We now assume that $r \leq 1$.

As in \cite{stevens}, let 
$$\Sigma_0(p) = \{ \left(\begin{matrix}a & b \\ c & d \end{matrix}\right) \in M_2(\Z_p),\  p \not | a,\  p |c,\  ad-bc\neq 0 \}.$$ 
This is a monoid that contains the standard Iwahori group at $p$ as a submonoid.
If $w \in \WW(\Q_p)$ is a point on the weight space, and
$\tw : \Z_p^\ast \rightarrow \Q_p^\ast$ is the corresponding character, then following Stevens (\cite{stevens}) 
we define the {\it weight-$w$ actions} of $\Sigma_0(p)$ on 
$\A[r]$ and $\D[r]$ by
\begin{eqnarray*} 
(\gamma \cdot_w f)(z) & = & \tw(a+cz) f\left(\frac{b+dz}{a+cz}\right),\ \ \ f \in \A[r] \\
(\mu \cdot_w \gamma)(f) & = & \mu (\gamma \cdot_w f),  \ \ \ \mu \in \D[r]
\end{eqnarray*}
where $\gamma=\left(\begin{matrix}a & b \\ c & d \end{matrix}\right)$. 
Some terms in the above definitions deserve explanations. Since $\tw$ is analytic on each of the
"components" $a+p \Z_p$ of $\Z_p^\ast$, and since $p | c$, $\tw(a+cz)$ may be seen as a power series in 
$\A=\A[1]$, that we may see also as an element of $\A[r]$ since $r \leq 1$. Also since $r\leq 1$, it is clear that $\frac{b+dz}{a+cz}$ belongs to $B[\Z_p,r]$ if $z$ does, so $f\left(\frac{b+dz}{a+cz}\right)$ also belongs to $\A[r]$. 

When we think of $\A[r]$ and $\D[r]$ as provided with the weight-$w$ action of $\Sigma_0(p)$, we denote them by $\A_w[r]$ and $\D_w[r]$. 
Note that when $w = k$ is integral, the weight-$k$ actions are exactly the ones defined in the introduction (\S\ref{intmodsymb}).
\par \bigskip

We shall now define an action of $\Sigma_0(p)$ on the $R$-modules $R \hotimes_{\Q_p} \A[r]$ and $R \hotimes_{\Q_p} \D[r]$ that interpolate the various  weight-$w$ actions for $w \in W = \sp R$. 

Note first that we can see $K(a+c z)$ (for $a,c \in \Z_p$, $p \not | a$, $p|c$, and $z$ a formal indeterminate) as an element of $R \hotimes \A$.
Indeed, up to shrinking $R$ if necessary, we may assume that $R \simeq \Q_p\{T\}$
and $K : \Z_p^\ast \rightarrow R^\ast$ is given by 
$$K(x)=\teich{x}^\beta \left( \frac{x}{\teich{x}} \right)^T$$
for some $\beta=1,\dots,p-1$, 
where for $y \in 1+p\Z_p$,
$y^T := \exp(T \log y)$ and $\teich{x}$ is the Teichmuller representative of $x$ mod $p$.
Therefore, noting that $\teich{a+cz} = \teich {a}$ for $z \in \Z_p$, it is natural to set
$$K(a+cz) = \teich{a}^\beta \exp (T \log (\frac{a+cz}{\teich{a}}))$$
as a formal series in $T$ with coefficients in $\A[1]$ (since $\log (\frac{a+cz}{\teich{a}})$
 belongs to $p \A[1]$) and it is easily seen that this series defines an element 
 of $R \hotimes \A[1]$. If $r \leq 1$ we can also see $K(a+cz)$ as an element of $R \hotimes \A[r]$. 
 
Now, let $\gamma = \left(\begin{matrix}a & b \\ c & d \end{matrix}\right)$. We define a morphsim of $R$-modules
$\gamma:\ R \otimes \A[r] \rightarrow R \hotimes \A[r]$ by sending $1 \otimes f$ to $\gamma \cdot (1 \otimes f) = K(a+cz) 
\left(1 \otimes f\left(\frac{b+dz}{a+cz}\right) \right)$. It is clear that this morphism is continuous, and therefore extends
to an endomorphism $\gamma: R \hotimes \A[r] \rightarrow R \hotimes \A[r]$. This defines an action of $\Sigma_0(p)$ on the $R$-modules $R \hotimes \A[R]$, and also a dual action on $R \hotimes \D[r]$ (using Lemma~\ref{ON}(iii)).

\begin{lemma} \label{lemmaspec1}
If $w \in \sp R(\Q_p)=\Hom(R,\Q_p)$ then there are canonical isomorphisms of $\Q_p$-Banach spaces
\begin{eqnarray*} (R \hotimes \A[r]) \otimes_{R,w} \Q_p &\rightarrow& \A_{w}[r]\\ 
(R \hotimes \D[r]) \otimes_{R,w} \Q_p &\rightarrow& \D_{w}[r]\end{eqnarray*}
 compatible with the action of $\Sigma_0(p)$.
\end{lemma}
\begin{pf} That there is a natural isomorphism of $\Q_p$-Banach spaces as in the first line is obvious, and we need only
to check that the actions of $\Sigma_0(p)$ are compatible. What we need to check for the first isomorphism 
is that if $f  \in \A[r]$, the image by the map $w \hotimes \Id: R \hotimes \A[r] \rightarrow \A[r]$ of
$\gamma \cdot (1 \otimes f)$ is $\gamma \cdot_w f$, that is that it sends  $K(a+cz)\left(1 \otimes f\left(\frac{b+dz}{a+cz}\right) \right)$ to 
$\tw(a+cz) f\left(\frac{b+dz}{a+cz}\right)$.
This is clear using (\ref{Kw}). The second isomorphism follows from the first.
\end{pf}
 
\label{actions}

\subsection{Distributions-valued modular symbols and their families}
\label{dismodsymb}

\subsubsection{Reminder on modular symbols (after Stevens \cite{stevens})}
\label{remindermodsymb}

From now, we fix an integer $N$ prime to $p$ and set $\Gamma = \Gamma_1(N) \cup \Gamma_0(p) \subset \Gl_2(\Q)$.

Let $S_0(p)$ be the submonoid of $\Sigma_0(p)$ consisting of matrices whose entries are rational integers.
Note that $\Gamma \subset \Gamma_0(p) \subset S_0(p) \subset \Gl_2(\Q)$. 

Let $\Delta_0$ be the abelian group of divisors of degree $0$ on $\PP^1(\Q)$. It is provided with a natural action of $\Gl_2(\Q)$, hence by restriction, of $S_0(p)$, and in particular of $\Gamma$. 

For $A$ any abelian group with a right action of $S_0(p)$, the abelian group $\Hom(\Delta_0,A)$ has a natural action of $S_0(p)$ coming from the actions of $S_0(p)$ on $\Delta_0$ and $A$ (namely $\Phi_{|\gamma}(\delta)=\Phi(\gamma \cdot \delta)_{|\gamma}$ for $\Phi \in \Hom(\Delta_0,A)$ and all $\delta \in \Delta_0$.)  
Recall that one defines the $A$-valued modular symbols
with level $\Gamma$ as $$\Symb_\Gamma(A) := \Hom_\Gamma(D_0,A).$$
By its very definition as the abelian group of $\Gamma$-invariants in  the larger abelian group with $S_0(p)$-action $\Hom(D_0,A)$, 
the abelian group $\Symb_\Gamma(A)$ inherits an action of the Hecke algebra of $S_0(p)$ w.r.t $\Gamma$.
This algebra contains the classical Hecke operators $T_l$ for all prime $l$ not dividing $Np$, the \diamonds\  as well as the classical Atkin-Lehner 
operator $U_p$, that is to say it contains $\HH$. Since the matrix $\diag(-1,1)$ belongs to $S_0(p)$ and normalizes $\Gamma$, it also defines an involutive operator on $\Symb_\Gamma(A)$ commuting with the action of $\HH$. We denote by 
$\Symb_\Gamma^{\pm}(A)$ the subgroup of $\Symb_\Gamma(A)$ of elements that are fixed or multiplied by $-1$
by that involution. If $2$ acts invertibly on $A$, we have $\Symb_\Gamma(A) = \Symb_\Gamma^+(A) \oplus \Symb_\Gamma^-(A)$. 

Of course, if $A$ has a structure of module over any commutative ring $B$ such that the $S_0(p)$-action is $B$-linear, then $\Symb_\Gamma(A)$ (resp. $\Symb_\Gamma^\pm(A)$) has a natural structure of $B$-module and the Hecke operators defined above are $B$-linear 
applications. If $B$ is a Banach algebra and $A$ is a Banach $B$-module, such that $S_0(p)$ acts by
continuous endomorphisms and $\Gamma$ by unitary endomorphisms, then $\Symb_\Gamma(A)$ has a structure of Banach $B$-module defined by 
$\NN{\Phi} = \sup_{D \in \Delta_0} \NN{\Phi(D)}$  (the $\sup$ is well defined because $\NN{\Phi(D)}$ depends only on the class of $D$ on $\Delta_0/\Gamma$ by our assumption and because $\Delta_0/\Gamma$ is finite)
for $\Phi \in \Symb_\Gamma(A)$, on which 
the Hecke operators act continuously. 

 Let $\delta_i$, $i=1,\dots,h$ be representatives in $\Delta_0$ of the elements of 
 $\Delta_0/\Gamma$, 
and let $\Gamma_i$ be the stabilizer of $\delta_i$. Those groups are well-known to be finite.
Then we have a natural isomorphism of groups (resp. of $B$-modules, resp. of normed $B$-modules 
in the two situations considered the preceding paragraph) 
\begin{eqnarray}\label{symb} \Symb_\Gamma(A) \rightarrow \oplus_{i=1}^h A^{\Gamma_i}\end{eqnarray}
given by $\Phi \mapsto (\Phi(\delta_i))_{i=1,\dots,h}$.
This isomorphism makes clear that when $A$ is a $B$-module and $B$ has characteristic $0$, then the formation
of $\Symb_\Gamma(A)$ commutes with arbitrary base changes $B \rightarrow B'$.

\subsubsection{Families of distributions-valued modular symbols}

We keep the notations of the preceding \S, and we now consider the
$R$-modules $\Symb_\Gamma^\pm(R \hotimes \D[r])$. We summarize its properties in
\begin{prop} \label{propsymb}
 For any choices of sign $\pm$,
\begin{itemize}
 \item[(i)] the group $\Symb_\Gamma^\pm(R \hotimes \D[r])$ is a Banach $R$-module which 
 satisfies Buzzard's property $(Pr)$;
\item[(ii)] The elements of $\HH$ act continuously on this module, and $U_p$ acts compactly;
\item[(iii)] If $w \in W(\Q_p) = \Hom(R,\Q_p)$,
then there is a natural isomorphism of $\Q_p$-Banach spaces, compatible with the action of $\HH$
$$\Symb_\Gamma^\pm(R \hotimes \D[r])  \otimes_{R,w} \Q_p \simeq \Symb_\Gamma^\pm(\D_{w}[r]).$$
\end{itemize}
\end{prop}
\begin{pf} The point (i) is clear by (\ref{symb}) and Lemma~\ref{ON}(i). The continuity of the actions of elements of $\HH$ is clear by the discussion
on the topology above and the compactness of $U_p$ follows from the fact that it factors 
through $\Symb_\Gamma^\pm(R \hotimes \D[r/p])$ (see \cite{stevens}), combined
with Lemma~\ref{ON}(ii) and (\ref{symb}). The point (iii) is a special case, in view of Lemma~\ref{lemmaspec1}, of the compatibility of the 
formation of $\Symb_\Gamma(A)$ with base change shown above. 
\end{pf}

Since $U_p$ acts compactly on $\Symb_\Gamma^\pm(R \hotimes \D[r])$ which is (Pr), Riesz' theory (as extended by Serre, Coleman and Buzzard: see \cite[chapter 2]{buzzard})
applies to it, and allows us to define, for any positive real number $\nu$ a sub-module $\Symb_\Gamma^\pm(R \hotimes \D[r])^{\leq \nu}$ on which
$U_p$ acts with generalized eigenvalues of valuation $< \nu$. This submodule is finite and projective over $R$ by \cite[Lemma 2.11]{buzzard}.
It is stabilized by $\HH$ since $U_p$ commutes with $\HH$.

Similarly, for any $w \in W(\Q_p)$, $U_p$ acts compactly on $\Symb_\Gamma^\pm(\D_w[r])$ (either by a direct argument, or by (ii) and (iii) of the above Proposition), and we can consider the finite-dimensional $\Q_p$-space, with an $\HH$-module structure $\Symb_\Gamma^\pm(\D[r])^{\leq \nu}$.
By (iii) of the proposition above, and elementary reasoning with the Riesz decomposition (cf. \cite[\S4]{stevenspollack2}), we have
\begin{cor} \label{corsymb} 
If $w \in W(\Q_p) = \Hom(R,\Q_p)$,
then there is a natural isomorphism of $\Q_p$-vector spaces compatible with the action of $\HH$
$$\Symb_\Gamma^\pm(R \hotimes \D[r])^{< \nu}  \otimes_{R,w} \Q_p \simeq \Symb_\Gamma^\pm(\D_{w}[r])^{< \nu}.$$
 \end{cor} 

In what follows, we fix an $r$ such that $0 < r \leq 1$ ($r=1$ is the obvious choice, with $r=1/p$ as in the introduction a serious competitor)  
hence dropping the $[r]$ from the notations: $$\D:=\D[r].$$

\subsubsection{Stevens' control theorem and the Eichler-Shimura isomorphism}

\label{stevenscontrol}

Let $k \geq 0$ be an integer. Following Stevens (using notations of Stevens and Pollack), we call 
$V_k$ the set of homogeneous polynomial of degree $k$, and we see it as a right $\Gl_2(\Q_p)$-module by setting
$$(P|\gamma)(X,Y)=P(dX-cY,-bX+aY)\ \ \ \text{ if }\gamma=\left(\begin{matrix} a & b \\ c & d\end{matrix}\right).$$
There is an exact sequence (cf. \cite{stevens}) 
\begin{eqnarray*}  0 \longrightarrow \D_{-2-k}(k+1) \stackrel{\Theta_k} \longrightarrow 
\D_k \stackrel{\rho_k} \longrightarrow V_k  \longrightarrow 0 \end{eqnarray*}
of $\Q_p$-spaces with actions of $S_0(p)$ (the $(k+1)$ in $\D_{-2-k}(k+1)$ indicates a twist by the $k+1$-power of the determinant in that action), where 
$\rho_k^\ast(\mu) =  \int (Y-zX)^k\  d\mu(z)$
and $\Theta_k$ is the dual of the $k$-th derivative map,
which induces an exact sequence
\begin{eqnarray*}  0 \longrightarrow \Symb_\Gamma(\D_{-2-k})(k+1) 
\stackrel{\Theta_k^\ast} \longrightarrow 
\Symb_\Gamma(\D_k) \stackrel{\rho_k^\ast} \longrightarrow \Symb_\Gamma(V_k).
 \end{eqnarray*}
The morphisms $\Theta_k^\ast$ and  $\rho_k^\ast$ respect the action of $\HH$ and the involution $\iota$.

\begin{theorem}[Stevens' control theorem]
The morphism $\rho_k^\ast$ is surjective, and it induces an isomorphism (of $\Q_p$-spaces with an action of $\HH$ and $\iota$)
$$\rho_k^\ast: \Symb_\Gamma(\D_k)^{< k+1} \rightarrow \Symb_\Gamma(V_k)^{<k+1}$$
\end{theorem}
\begin{pf} See \cite[Theorem 7.1]{stevens}.
\end{pf}
This theorem is usually coupled (compare \cite{stevenspollack2}) with this form of  the Eichler-Shimura isomorphism:
\begin{prop}[Eichler-Shimura] \label{es} There exists a natural $\HH$-isomorphism
\begin{eqnarray*} \es: \Symb_\Gamma(V_k)& \isomo & M_k(\Gamma) \oplus 
S_k(\Gamma)
\end{eqnarray*}
It splits according to the $\iota$ involution as a sum $\es = \es^+ \oplus \es^-$ into two injective $\HH$-morphisms
\begin{eqnarray*}
\es^+: \Symb_\Gamma^+(V_k) &\hookrightarrow& M_k(\Gamma) \\
\es^-: \Symb_\Gamma^-(V_k) &\hookrightarrow& M_k(\Gamma)
\end{eqnarray*}
whose images contains both $S_k(\Gamma)$.

The inverse image by $\es$ of the space of Eisenstein series in $M_k(\Gamma)$ is
the space of {\rm boundary modular symbols}, that is modular symbols that extend
to  $\Gamma$-invariant morphisms from $\Delta$ to $V_k$, where $\Delta$ is the set of all divisors on $\PP^1(\Q)$, not necessarily of degree $0$.
\end{prop}
\begin{pf} This theorem is quoted in \cite{merel},
\cite{stevens}, \cite{stevenspollack2}, but I have been unable to find a really
complete proof in the literature. We shall admit it.

The cuspidal part of the statement is due to Shokurov \cite{shokurov}
(generalizing results of Manin in the level $1$ case), and details can be found in \cite{merel} and \cite{steinthesis}. Alternatively, following an argument due to Ash and Stevens, this cuspidal case can be obtained
by combining the classical Eichler-Shimura isomorphisms (involving $H^1_c(Y(\Gamma),\tilde V_k)$, cf. e.g. \cite[\S6.2]{hida}) 
with isomorphisms coming from exact sequences in cohomology.
Indications relative to how to proceed with the Eisenstein case are given in \cite{merel}.

For later use, we recall how are defined the inverse of the isomorphism $e_s^\pm$ on cuspidal modular form: for a cuspidal eigenform $f$, $(\es^\pm)^{-1}(f)$ is the $\pm$-part of the
$V_k$-valued modular symbol that to the divisor $\{r\} -\{r'\} \in \Delta_0$ attaches the values $\int_r^{r'} 2 \pi i f(z) (X-zY)^k \ dz$ divided by the period $\Omega_f^\pm$.
\end{pf} 
\begin{cor}\label{corcontrol} There exist canonical $\HH$-morphism 
\begin{eqnarray*}
\es^+ \circ  \rho_k^\ast : \Symb_\Gamma^+(\D_k)^{< k+1} & \hookrightarrow & M_{k+2}(\Gamma)^{<k+1} \\
\es^- \circ \rho_k^\ast : \Symb_\Gamma^-(\D_k)^{< k+1} & \hookrightarrow & M_{k+2}(\Gamma)^{< k+1}
\end{eqnarray*}
whose images contain $S_{k+2}^{< k+1}(\Gamma)$.
\end{cor}

\subsubsection{The sign of an Eisenstein series}

\label{signeisenstein}

Let $f$ be a new Eisenstein series of level $\Gamma_1(N)$, and $f_\beta$ one of its refinement. We shall say that the sign of $f_\beta$ is $+1$ (resp. $-1$) if $f_\beta$ belongs to the image of $\es^+$ (resp. of $\es^-$). By Prop.~\ref{es}, the sign of 
$f_\beta$ is well-defined, and we shall denote if $\epsilon(f_\beta)$. 

\begin{remark} It is easy to see that Eisenstein series of level $\Gamma_0(N)$ 
with $N$ square free, have $\epsilon(f_\beta)=\epsilon(f_\beta)=1$. But in general,
even for $\Gamma_0(N)$ with $N=9$ and $p=2$, there are forms $f_\alpha$ with $\epsilon(f_\alpha)=-1$. 
\end{remark}

\subsection{Mellin transform of distribution-valued modular symbols}

This  won't be needed until \S\ref{sectionL}.

\subsubsection{Overconvergent distribution and Mellin transform in families}

\label{mellin}
This section will not be needed until \S\ref{sectionL}.

Let us call $\D[0]=\projlim \D[r]$, which is the Frechet space of so-called {\it overconvergent distributions}.
For every $w \in W(\Q_p)$, it has a natural weight-$w$ action of $\Sigma_0(p)$, coming from the 
compatible weight $w$-action on $\D[r]$. We call $\D[0]_w$ the space $\D[0]$ endowed with that action.
Similarly $R \hotimes \D[0]$ has a natural action of $\Sigma_0(p)$ coming from the actions on $R \hotimes \D[r]$
defined in~\S\ref{actions}.

\begin{definition} We call $\RR$ the $\Q_p$-algebra of analytic functions on $\WW$.
\end{definition}
 The space $\RR$ is naturally a Frechet space. Elements of $\RR$ can be seen as $\C_p$-valued functions of one variable $\sigma$ belonging in $\WW(\C_p)$.
\begin{definition} \label{defmellin} 
The {\it Mellin transform} is the application $$\Mel:  \D[0] \rightarrow \RR$$ that sends a distribution $\mu$ to
the function $\sigma \mapsto \int_{\Z_p^\ast} \tilde \sigma d\mu$ for $\sigma \in \WW(\C_p)$, $\tilde \sigma$ being as usual the corresponding character $\Z_p^\ast \rightarrow \C_p^\ast$. The {\it Mellin transform in family} is the application $\Mel_R: R \hotimes \D[0] \rightarrow R \hotimes \RR$ 
defined as $\Id_R \hotimes \Mel$.
\end{definition} 
 Note that the Mellin transforms $\Mel$ is  continuous for the Fr\'echet topologies. This is used to see that the definition of $\Mel_R$ 
makes sense.

The elements of $R \hotimes \RR$ can be seen as  two-variables analytic functions on $W \times \WW$. The first variable, in $W$, is usually denoted by $w$ or $k$, and the second variable by  $\sigma$ or $s$. We have the following obvious (from the properties of completed tensor products) 
compatibility relation: If $\mu \in R \hotimes \D[0]$, and if $w \in W(\Q_p)=\Hom(R,\Q_p)$, we denote by $\mu_w$ the image of $\mu$ by the 
map $w \otimes \Id_{\D[0]} : R \hotimes \D[0] \rightarrow \D[0]_w$. Then for all $\sigma \in W(\C_p)$, we have
\begin{eqnarray} \label{compMel} \Mel_R(\mu)(w,\sigma) = \Mel(\mu_w)(\sigma) \end{eqnarray}

\subsubsection{From distribution-valued modular symbols to analytic functions}
\label{from}
\label{lambda}

Fix $\nu >0$. In this \S, we shall use the Mellin transforms to define two natural maps
$$\Lambda_R: \Symb_\Gamma(R \hotimes \D)^{\leq \nu} \rightarrow  R \hotimes \RR$$
and for every $w \in W(\Q_p)$, a map
$$\Lambda: \Symb_\Gamma(\D_w)^{\leq \nu} \rightarrow \RR.$$
The two maps are compatible in the sense that for $\phi \in  \Symb_\Gamma(R \hotimes \D)^{\leq \nu}$, if $\phi_w$ denotes its image in 
$\Symb_\Gamma(\D_w)^{\leq \nu}$, we have for all $\sigma \in \WW(\C_p)$,
\begin{eqnarray} \label{complambda} \Lambda_R(\phi)(w,\sigma)= \Lambda(\phi_w)(\sigma). \end{eqnarray}

We  follow closely a construction of Stevens.

According to \cite[Prop. 5.6]{stevens} the natural map 
\begin{eqnarray} \label{stev1} \Symb_\Gamma(\D[0]_w)^{\#} \isomo \Symb_\Gamma(\D_w)^{\#} \end{eqnarray}
induced by the restriction map $\D[0]_w \rightarrow \D_w$ is an isomorphism.
Here, for $V$ any vector spaces on which $U_p$ acts, $V^\#:=\cap_{n=1}^\infty \U_p^n(V)$. Therefore, if $V$ is
a Banach space satisfying $(Pr)$, $V^\# = \cup_{\nu>0} V^{\leq \nu}$.

Exactly the same proof as in \cite{stevens} gives an analog isomorphism in families:
\begin{eqnarray} \label{over} \Symb_\Gamma(R \hotimes \D[0])^{\#} \isomo  \Symb_\Gamma(R \hotimes \D )^{\#} \end{eqnarray}
which is compatible with the $\HH$-action and with~(\ref{stev1}) via the specialization maps $\D \hotimes R \rightarrow \D_w$
and $\D[0] \hotimes R \rightarrow \D[0]_w$.

If $\phi \in \Symb_\Gamma(R \hotimes \D)^{\leq \mu}$, then $\phi \in \Symb_\Gamma(R \hotimes \D)^\#$, so is identified with an
 element $\tilde \phi$ of 
$\Symb_\Gamma(R \hotimes \D[0])^\#$ by isomorphism~(\ref{over}). The evaluation of $\tilde \phi$ at the divisor $\{\infty\}-\{0\}$
gives an element $\tilde \phi(\{\infty\}-\{0\})$ in $R \hotimes \D[0]$ of which we can take the image by the application $\Mel_R$.
We set $$\Lambda_R(\phi):=\Mel_R(\tilde \phi(\{\infty\}-\{0\})).$$ The construction of $\Lambda$ is similar (and slightly simpler): 
$$\Lambda(\phi):=\Mel(\tilde \phi(\{\infty\}-\{0\})).$$ The compatibility~(\ref{complambda}) is obvious from the definitions and (\ref{compMel}).

\subsubsection{Properties of $\Lambda$}

Let us recall the following well-known result:
\begin{prop}[Visik, Stevens] \label{proporder} If $\phi \in \Symb_\Gamma(\D_w)^{\leq \rho}$, then $\Lambda(\phi)$ is, as a function of $\sigma$, of order less that $v_p(\rho)$.
\end{prop}
The definition of the {\it order} is recalled in the introduction.

We now turn to interpolation properties. 
\begin{lemma} \label{interpolation} Let $k$ be a non-negative integer. Let $\phi \in \Symb_\Gamma(\D_k)$ be an eigensymbol
which is a generalized eigenvector for $U_p$, of generalized eigenvalue $\beta \neq 0$. We assume that $\rho_k^\ast((\U_p-\beta)(\phi)) = 0$.
Let $\tilde \phi$ be the  modular symbol in $\Symb_\Gamma(\D[0]_k)^{\leq \nu}$ corresponding to $\phi$, 
and let $\mu = \tilde \phi (\{\infty\}-\{0\}) \in \D[0]$.
 Then for every integer such that $0 \leq j \leq k$, every integer $n \geq 0$, and every integer $a$ such that $0 < a < p^n$ we have
\begin{eqnarray} \label{calcul} \int z^j {\bf{1}}_{a+p^n\Z_p} \ d \mu(z) = \beta^{-n} \phi(\{\infty\}-\{a/p^n\})((p^n z+a)^j).\end{eqnarray}
If we write $\rho_{k}^\ast(\phi)(\{\infty\}-\{a/p^n\})=\sum_{i=0}^k c_{i} X^i Y^{k-i} \in V_k$,
then we have 
\begin{eqnarray} \label{calculbis} \int z^j {\bf{1}}_{a+p^n\Z_p} \ d \mu(z) =
\beta^{-n} \sum_{i=0}^j c_{i} (-1)^{k-i} \frac{j!(k-i)!}{k!(j-i)!} p^{ni}a^{j-i}
\end{eqnarray}

\end{lemma}
\begin{pf} 
The equation (\ref{calculbis}) follows from (\ref{calcul}) by the definition of $\rho_k^\ast$ and a simple computation.

Let $r$ be the smallest integer such that $(U_p-\beta)^r \phi = 0$. We shall prove 
(\ref{calcul}) by induction on $r$.

If $r=1$, $\phi$ is a true eigenvector for $U_p$, and 
it is a relatively simple computation, done in \cite[\S6.3]{stevenspollack1} that
\begin{eqnarray} \label{calc2} \int z^j {\bf{1}}_{a+p^n\Z_p}\ d\mu(z) = \beta^{-n} \phi(\{\infty\}-\{a/p^n\})((p^nz+a)^j).\end{eqnarray}

Let us assume that $r>1$, and that (\ref{calcul}) (and hence (\ref{calculbis})) are known for all $r'<r$. For any $n$ we have $(U_p^n-\beta^n)\phi=(U_p-\beta)Q(U_p)\phi$ where $Q$ 
is some obvious polynomial. Set $\phi'=(U_p-\beta)Q(U_p)\phi$. Then obviously $(U_p-\beta)^{r-1}\phi'=0$ and by assumption $\rho_k^\ast(\phi')=0$, which
implies{\it  a fortiori} $\rho_k^\ast ((U_p-\beta)\phi') = 0$, so we can apply the induction hypothesis to $\phi'$, which gives us, calling $\mu'$ the distribution $\tilde \phi'(\{\infty\}-\{0\})$:
\begin{eqnarray} \label{calc3} \int z^j {\bf{1}}_{a+p^n\Z_p}\ d\mu'(z) = 0.\end{eqnarray}
But independently the same computation as in \cite[\S6.3]{stevenspollack1}, starting with $U_p^n \phi - \phi' = \beta^n \phi$ give us
that 
$$\int z^j {\bf{1}}_{a+p^n\Z_p}\ d\mu(z)  - \beta^{-n} \int z^j {\bf{1}}_{a+p^n\Z_p}\ d\mu'(z) = \beta^{-n} \phi(\{\infty\}-\{a/p^n\})((p^nz+a)^j).$$
Adding  $\beta^{-n}$ times (\ref{calc3}) we get (\ref{calc2}) and we conclude as above.
\end{pf}

\begin{prop} \label{inter0} Let $k$ be a non-negative integer. Let $\phi \in \Symb_\Gamma(\D_k)$ be an eigensymbol
which is a generalized eigenvector for $U_p$, of generalized eigenvalue $\beta \neq 0$. We assume that $\rho_k^\ast(\phi) = 0$.
Then for all character $\varphi : \Z_p^\ast \rightarrow \C_p^\ast$ of finite image and all integers $j$ with $0 \leq j \leq k$, we have
$$\Lambda(\phi)(\varphi t^j)=0$$
\end{prop}
\begin{pf} With our assumptions the above lemma tells us that for all integers $j$, $n$ and $a$ $0 \leq j \leq k$, $n \geq 0$ an integer,  and integer $a$ such that $0 < a < p^n$ we have $\int z^j {\bf{1}}_{a+p^n\Z_p} \ d \mu(z) = 0$. But $\Lambda(\phi)(\varphi t^j)=\int \psi(z) z^j d \mu(z)$
is defined  as a limit of Riemann sums that are linear combinations on terms of the form $\int z^j {\bf{1}}_{a+p^n\Z_p} \ d \mu(z)$ with $j,n,a$ as above. So the Riemann sums, and their limit, are all $0$.
\end{pf}

\begin{prop} \label{interf} Let $k$ be a nonnegative integer, $f_\beta$ the refinement of a cuspidal newform in $S_{k+2}(\Gamma)$. 
Let $\phi$ be a non-zero symbol in $\Symb_\Gamma^\pm(\D_k)_{(f_\beta)}$ such that $\rho_k^\ast(\phi) \neq 0$.
Then for any special character of the form $\sigma=\varphi t^j$, with $0 \leq j \leq k$, we have
$$\Lambda(\phi)(\sigma) = e_p(\beta,\psi t^j) \frac{m^{j+1}}{(-2\pi i)^j} \frac{j!}{\tau(\varphi^{-1})} \frac{L_\infty(f\varphi^{-1},j+1)}{\Omega_f^{\pm}},$$ up to a non-zero scalar depending only on $\phi$, not on $\sigma$, with $\L_\infty$, $\Omega_f^{\pm}$ $e_p(\beta,\varphi t^j)$, $m$ as in the introduction.
\end{prop}
\begin{pf} 
Since $\rho_k^\ast(\phi)\neq 0$, and $\es^\pm$ is an isomorphism, $\es^\pm \circ \rho_k^\ast(\phi)$ is a non-zero vector in $M_\Gamma(k+2)_{(f_\beta)}$. Since this space
has dimension $1$ (Lemma~\ref{lemmanewform}(i)), we see that $\es^\pm \circ \rho_k^\ast(\phi)=f_\beta$ up to a non-zero scalar.
Also, $(U_p -\beta)f_\beta=0$ so we see that $\rho_k^\ast((U_p-\beta)\phi)=0$, which allows us to apply Lemma~\ref{interpolation}.
The lemma then follows by a standard computation from~(\ref{calcul}) combined with the explicit description of $\es^\pm$ for a cuspidal form $f_\beta$ recalled in the proof of Proposition~\ref{es}.
\end{pf}

\subsection{Local pieces of the eigencurve}

\subsubsection{Construction of the local pieces}

\label{localpieces}

We keep the notations of the preceding paragraphs, in particular $W =\sp R$ is an open affinoid subspace of the weight space $\WW$, and $\Gamma = \Gamma_0(p) \cap \Gamma(N)$ with $N$ an integer prime to $p$. We choose a positive real number $\nu$. Our aim is to construct by means of distributions-valued modular symbols an affinoid
space $\CC^\pm_{W,\nu}$ together with a finite flat ``weight'' morphism 
$\kappa: \CC^\pm_{X,\nu} \rightarrow W$.

The construction follows the general method described in \cite{buzzard}:
we consider the space $\Symb_\Gamma^\pm(R \hotimes \D)^{\leq \nu}$ introduced in \S\ref{famdis}
This is by definition a finite projective $R$-module on which $\HH$ acts.
\begin{definition} We call $\T_{W,\nu}^\pm$ the $R$-subalgebra of 
$\End_R(\Symb_\Gamma^\pm(R \hotimes \D)^{\leq \nu})$ generated by the image of $\HH$.
\end{definition}  The algebra $\T_{W.\nu}^\pm$ is obviously finite and torsion-free as an $R$-module, and since $R$ is Dedekind, 
it is also flat. As a finite algebra over $R$, it is an affinoid algebra, and we set $$\CC^{\pm}_{W,\nu}:= \sp \T^\pm_{W,\nu}.$$ 
The weight morphism
$$\kappa^\pm: \CC_{W,\nu}^\pm \rightarrow W=\sp R$$
 is simply the map corresponding to the structural morphism $R \rightarrow \TT_{W,\nu}^\pm$ and is therefore finite and flat. This completes the construction.

\subsubsection{The ``punctual eigencurve''}

Let $w \in W(\Q_p) = \Hom(R,\Q_p)$ be some point on the weight space. We can perform the same construction
in smaller at the point $w$: 
\begin{definition} We define $\T_{w,\nu}^\pm$ as the $\Q_p$-subalgebra 
of $\End_{\Q_p}(\Symb_\Gamma^\pm(\D_w)^{\leq \nu})$ generated by the image of $\HH$.
\end{definition}
By construction, there is a natural bijection between the $\bar \Q_p$-points of 
$\spec \T_{w,\nu}^\pm$ (that is, the morphisms of algebras $T_{w,\nu}^\pm \rightarrow \bar \Q_p$) and the systems of $\HH$-eigenvalues appearing
in $\Symb_\Gamma^\pm (\D_w)^{\leq \nu}$ (see \S\ref{notations}). If $x$ is such a point, or systems of $\HH$-eigenvalues, we denote by
$\Symb_\Gamma^\pm(\D_w)_{(x)}$ the corresponding {\bf generalized eigenspace} (we drop the ${}^{\leq \nu}$ which is redundant since
being in the $(x)$-generalized eigenspace implies being a generalized eigenvector for $U_p$ of generalized eigenvalue $x(U_p)$ which by definition has valuation $\leq \nu$). Similarly we note $(\T_{w,\nu}^\pm)_{(x)}$ the localization of $\T_{w,\nu}^\pm \otimes_{\Q_p} \bar \Q_p$ at the maximal ideal corresponding to $x$. Obviously we have $$\Symb_\Gamma^\pm(\D_w)_{(x)} = \Symb_\Gamma^\pm(\D_w)^{\leq \nu} \otimes_{T_{w,\nu}^\pm} (\T_{w,\nu}^\pm)_{(x)}.$$

If $w$ is an integral weight $k \in \N$ and  if the system of eigenvalues $x$ happens to be the system of eigenvalue of a refinement $f_\beta$
of a  modular newform $f$ of level $\Gamma_1(N)$ and weight (necessarily) $k+2$, recall (cf~\ref{notations}) that we sometimes write  $\Symb_\Gamma^\pm(\D_w)_{(f_\beta)}$
instead of  $\Symb_\Gamma^\pm(\D_w)_{(x)}$ and $(\T_{w,\nu}^\pm)_{(f_\beta)}$ instead of $(\T_{w,\nu}^\pm)_{(x)}$.

\subsubsection{Comparison global-punctual}

Let again $w \in W(\Q_p) = \Hom(R,\Q_p)$ be some point on the weight space.
Since by Prop.~\ref{propsymb}(iii), $\Symb_\Gamma^\pm(R \hotimes \D)  \otimes_{R,w} \Q_p \simeq 
\Symb_\Gamma^\pm(\D_{w})$, we have a natural {\it specialization} morphism of finite $\Q_p$-algebras
$$s_w : \T^\pm_{W,\nu} \otimes_{R,w} \Q_p \rightarrow \T^\pm_{w,\nu}$$
compatible with the morphisms from $\HH$.
\begin{lemma} \label{lemmasw} For all $w \in \W(\Q_p)$, $s_w$ is surjective, and its kernel is nilpotent.
In particular, $s_w$ induces a bijection between the $\bar \Q_p$-points of the fiber of $\kappa^\pm$ at $w$ and
of $\sp \T_{w,\nu}^\pm$:
$$(\kappa^\pm)^{-1}(w) (\bar \Q_p) \simeq \sp \T_{w,\nu}^\pm(\bar \Q_p).$$
\end{lemma}
\begin{pf} 
The surjectivity is obvious since the image of $s_w$ is a subalgebra of $\T_{w,\nu}$ containing the image of $\HH$, so is $\T_{w,\nu}$ by definition.
The assertion on the kernel is \cite[Lemme 6.6]{C1}. The last sentence is clear.
\end{pf}

Now if $x$ is a $\bar \Q_p$-point in $\spec \T_{w,\nu}^\pm$, it corresponds to a $\bar \Q_p$-point also denoted $x$ in $\CC_{W,\nu}^\pm$ such that $\kappa(x)=w$, and we can consider the rigid analytic localization $(\TT_{W,\nu}^\pm)_{(x)}$ of $\TT_{W,\nu}^\pm \otimes_{\Q_p} \bar \Q_p$ at the maximal ideal corresponding to $x$. This is naturally an algebra over $R_{(w)}$ defined as the rigid analytic localization of 
$R \otimes_{\Q_p} \bar\Q_p$ at the maximal ideal corresponding to $w$.

Localizing at  $x$ the map $s_w$ gives rise to a
 a surjective local morphism of finite local $\bar \Q_p$-algebra with nilpotent kernel
\begin{eqnarray} \label{sx} s_x : (\T_{W,\nu}^\pm)_{(x)} \otimes_{R_{(w)},w} \bar \Q_p \rightarrow  
(\T_{w,\nu}^\pm)_{(x)}
\end{eqnarray}

\subsubsection{Comparison with the Coleman-Mazur-Buzzard's eigencurve}

\begin{prop} \label{compCMB} There exist 
canonical closed immersions $$ \CC_{W,\mu} \hookrightarrow \CC^\pm_{W,\nu} \hookrightarrow \CC_{W,\nu}$$ which is compatible with the weight morphism $\kappa^\pm$ and $\kappa$ to the weight space $W$. For $k \in W(\Q_p)$ an integer,
 there exists $\HH$-injection between the $\HH$-modules 
 $$S^\dag_{k+2}(\Gamma)^{\ses, \leq \nu} \subset \Symb_\Gamma^\pm(\D_k )^{\ses, \leq \nu} \subset M^\dag_{k+2}(\Gamma)^{\ses, \leq\nu},$$ where $\ses$ means semi-simplification as an $\HH$-module.
 In particular, if $f_\beta$ is a cuspidal refined form, we have $\dim (S^\dag_{k+2}(\Gamma))_{(f_\beta)}= \dim (\Symb_\Gamma^\dag(\D_k))_{(f_\beta)}$.
 \end{prop}
 \begin{pf} This follows from \cite[Th\'eor\`eme 1]{chenevierJL} using both Coleman's control theorems (\cite{coleman} and Corollary \S\ref{controlcusp}) and
and Stevens's control theorems (cf. Corollary \S\ref{corcontrol}) to compare te obvious
 {\it classical structure} in the sense of (\cite{chenevierJL}) of $\CC_{W,\nu}$, $\CC^0_{W,\nu}$ and $\CC^\pm_{W,\nu}$ respectively.
 \end{pf}
\begin{remark} We note that the last assertion is a slightly 
more precise version of an unpublished theorem of Stevens (\cite{stevensfam}, quoted in \cite[Theorem 7.1]{stevenspollack2}) which states that the systems
of cuspidal $\HH$-eigenvalues appearing in $\Symb_\Gamma^\pm(\D_k)$ and $S^\dag_{k+2}(\Gamma)$ are the same. 
This method of proving such a result was suggested to me by G. Chenevier. 
\end{remark}

\section{$p$-adic L-functions}

\label{sectionL}

We keep the notations of \S\ref{localpieces}.

\subsection{The space $\Symb_\Gamma^\pm(R \hotimes \D)^{\leq \nu}$ as a module over 
$\T_{W,\nu}^\pm$}

We first prove an elementary lemma of commutative algebra.
\begin{lemma} \label{commalg} Let $R$ be a discrete valuation domain, $T$ be a finite free $R$-algebra,
and $M$ be  a finitely generated $T$-module,  free as an $R$-module. Suppose that $T$ is also a discrete valuation domain.
Then $M$ is free as a $T$-module as well.
\end{lemma}
\begin{pf} Since $T$ is a discrete valuation domain, we can write as a $T$-module $M=M_\free \oplus M_\tors$, with $M_\free$ a finite 
free module over $T$ of some rank $r$ and $M_\tors$ a finite torsion module over $T$. We claim that $M_\tors$ is torsion as an $R$-module: 
for any $m \in M_\tors$, there is an element $0 \neq x \in T$ such that $xm=0$; 
if $a_0$ is the constant term of a polynomial in $R[X]$ killing $x$, chosen of minimal 
degree, then $a_0 \neq 0$, and one sees easily that $a_0 m =0$, which proves the claim.
But since $M$ is locally free as an $R$-module, it has no torsion, so $M_\tors=0$, and 
$M=M_\free$ is free over $T$.
\end{pf}
\begin{remark} The hypothesis that $T$ is a discrete valuation domain cannot be weakened to ``$T$ is local'': 
let $R=k[[u]]$ where $k$ is a field of characteristic $0$, $M=R^2$, and $T$ be the subalgebra of $\End_R(M)$ generated over $R$ by the endomorphism $t$ of matrix (in the canonical basis) $\left( \begin{matrix} 0 & u \\ u^2 & 0 \end{matrix} \right)$. Then $t^2 = u^3$, and one sees that 
$T = R[t]/(t^2-u^3)$ as an $R$-algebra, so that $T$ is a finite free $R$ module of rank $2$, and as a ring is obviously local, but not regular. 
However $M$ is not free as a $T$-module, for if it were, it would be obviously of rank one, but at the maximal ideal $(t,u)$ of $T$, the fiber of $M$ is $M/(t,u)M = k^2$ has dimension $2$.

However, it can happen that $M$ is  free over $T$ without $T$ being 
a discrete valuation domain. For an example, just modify the preceding example by redefining $t$ as 
the endomorphism of matrix  $\left( \begin{matrix} 0 & 1 \\ u^3 & 0 \end{matrix} \right)$.
Then $T$ is the same algebra as before, but now $M$ is free over $T$, obviously generated by the second vector of the standard basis.
\end{remark}

\begin{prop} \label{proplocfree} Let $x$ be any smooth $\bar \Q_p$-point in $\CC_{W,\nu}^\pm$. Then 
the $(\TT_{W,\nu}^\pm)_{(x)}$-module $\Symb_\Gamma^\pm(R \hotimes \D)_{(x)}$ is free of finite rank.
\end{prop}
\begin{pf}
Since $x$ is smooth, and $\CC$ is a reduced curve, the local ring $T:= (\TT_{W,\nu}^\pm)_{x}$ is a discrete valuation domain.
Setting  $$M:=\Symb_\Gamma^\pm(R \hotimes \D) \otimes_{\T_{W,\nu}} T=\Symb_\Gamma^\pm(R \hotimes \D)_{(x)}$$ and replacing $R$ by its localization at $\kappa^\pm(x)$ 
bring us to the exact situation and notations of Lemma~\ref{commalg}. This lemma tells us that $M$ is free as a $T$-module.
\end{pf}

\begin{cor} \label{corsx} For any $w \in W(\Q_p)$, and $x \in \CC_{W,\nu}^\pm(\bar \Q_p)$ a smooth point 
such that $\kappa(x)=w$,  the morphism of $\bar \Q_p$-algebras
$$s_x : (\T_{W,\nu}^\pm)_{(x)} \otimes_{R_{(w)},w} \bar \Q_p \rightarrow  
(\T_{w,\nu}^\pm)_{(x)}$$
is an isomorphism.
\end{cor}
\begin{pf} 
Any operator $t \in (\T_{W,\nu}^{\pm})_{(x)}$ that acts
trivially on  $\Symb_\Gamma^\pm (R \otimes \D)_{(x)} \otimes_{R,w} \bar \Q_p = 
\Symb_\Gamma^\pm (\D_{w})_{(x)}$ must have $0$ image in $(\T_{W,\nu}^{\pm})_{(x)} \otimes_{R,w} \bar \Q_p$ since $\Symb_\Gamma^\pm (R \otimes \D)_{(x)}$ is free.
In other words, the map $s_x$ is injective. Since it is surjective by Lemma~\ref{lemmasw} and the discussion below it, it is an isomorphism.
\end{pf}

\begin{prop} \label{freerankone}
Let $x$ be a point in $\CC^\pm_{W,\nu}$. We assume that there exists a decent 
 newform $f$ of weight $k+2$ and level $\Gamma_1(N)$ and a refinement $f_\beta$ of $f$, such that the systems of eigenvalues of $f_\beta$ for $\HH$ corresponds to $x$. Then $\Symb_\Gamma^{\pm}(R \hotimes \D)_{(x)}$ is free of rank one on
$(\TT_{W,\nu}^\pm)_{(x)}$.
\end{prop}

\begin{pf} By Theorem~\ref{smooth} and Prop~\ref{compCMB}, $\CC^\pm_{W,\nu}$ is smooth at $x$. Therefore, on a suitable neighborhood of $x$ in $\CC^\pm_{W,\nu}$, all points are smooth, the map 
$\kappa$ is \'etale at all points but perhaps $x$ (by the openness of \'etaleness and the fact that $\CC^\pm_{W,\nu}$ is equidimensional of dimension $1$) and all
classical points corresponds to newform (by Lemma~\ref{newnew}).  
Since the localizations of a locally free finitely generated module have locally constant rank, we can evaluate
the rank $r$ of the free module $\Symb_\Gamma^{\pm}(R \hotimes \D)_{(x)}$  on
$(\TT_{W,\nu}^\pm)_{(x)}$ by replacing $x$ by any point in such a neighborhood. Choosing a classical point 
$y$ different from $x$, of weight $k'$
, we have that $r$ is the rank of $\Symb_\Gamma^{\pm}(R \hotimes \D)_{(y)}$ over 
$(\TT_{W,\nu}^\pm)_{(y)}$, and since $\kappa$ is \'etale at $y$, $r$ is the dimension of 
 $\Symb_\Gamma^{\pm}(\D_{k'})_{(y)}$ over $\bar \Q_p$. By Prop.~\ref{compCMB}
 this is the dimension of $M^\dag_{k'+2}(\Gamma)_{(y)}$, which is $1$ by Lemma~\ref{lemmanewform}(ii).
 \end{pf}

\subsection{Critical $p$-adic $L$-functions}

Let $f$ be a modular newform of weight $k+2 \geq 2$, level $\Gamma_1(N)$, Let $f_\beta$
be one of its two refinement, so $f_\beta$ is a modular form for $\Gamma=\Gamma_1(N) \cap\Gamma_0(p)$. Assume that $f_\beta$ is decent, and that $v_p(\beta) > 0$ in case $f$ is Eisenstein.

Let us fix a real number $\nu > k+1 \geq v_p(\beta)$, and as above an open affinoid subspace $W = \sp R$ of $\WW$ that contains $k$.
Let $\pm$ denotes any choice of the sign $+$ or $-$. 
Then we know that $f_\beta$ corresponds to a unique point $x \in \CC^0 \subset \CC$ (Lemma~\ref{lemmanewform}(iii)), and that this point
 is smooth (Theorem~\ref{smooth}). By  Prop~\ref{compCMB}, we see $x$ as a point of $\CC^\pm$ and this point $x$ actually belongs to the open set $\CC^\pm_{W,\nu}$.
We have $\kappa^\pm(x)=k$. 

Recall that $R_{(k)}$ is the rigid-analytic localization of $R \otimes_{\Q_P} \bar \Q_p$ at the maximal ideal corresponding to $k$. The ring $(\TT^\pm_{W,\nu})_{(x)}$ is 
an $R_{(k)}$-algebra. We call $e$ the index of ramification of $\kappa^\pm$ at $x$ (it is independent of $\pm$ by Corollary~\ref{corsmooth}), 
that is the length of the connected component containing $x$ of the fiber of $\kappa^\pm$ at $x$. In particular, $\kappa^\pm$ is \'etale at $x$ if and and only if $e=1$.

\subsubsection{Proof of Theorem~\ref{eigenspace}}

\label{proof1}

\begin{prop}\label{struc} For some choice of uniformizer $u \in R_{(k)}$
there is an isomorphism of $R_{(k)}$-algebras $R_{(k)}[t]/(t^e-u) \rightarrow (\TT^\pm_{W,\nu})_{(x)}$ sending $t$ to a uniformizer of  
$(\TT^\pm_{W,\nu})_{(x)}$.
\end{prop}
\begin{pf}  The rings $R_{(k)}$ and $(\TT^\pm_{W,\nu})_{(x)}$ are Henselian discrete valuation domains since they are the rigid analytic 
local rings of  rigid analytic curves at smooth points.
The algebra $(\TT_{W,\nu})_{(x)}$ is finite and flat (hence free) over $R_{(k)}$, since the morphism $\kappa$ is. We are therefore
dealing with a standard extensions of DVR, which is totally ramified of degree $e$, and moreover tamely ramified since
we are in characteristic $0$.  It is a well-known result (see e.g. 
\cite[Chapter 2, \S5]{lang}) that 
$(\TT^\pm_{W,\nu})_{(x)}$ 
has the form $R_{(k)}[t]/(t^e-u)$ for uniformizers $u$ and $t$.
\end{pf}

\begin{theorem}\label{thmL} The $\bar \Q_p$-space $\Symb_\Gamma^\pm(\D_k)_{(f_\beta)}$ has dimension $e$ and is a free module of rank one over
the algebra $(\TT^\pm_{k,\nu})_{(f_\beta)}$. This algebra  is isomorphic to $\bar \Q_p[t]/(t^e)$.
\end{theorem} 
\begin{pf} 
We have \begin{eqnarray*} 
(\TT^\pm_{k,\nu})_{(f_\beta)} & =& (\TT^\pm_{k,\nu})_{(x)}\text{\ \ by definition} 
\\ &=& (\TT^\pm_{W,\nu})_{(x)} \otimes_{R_{(k)},k} \bar \Q_p \text{\ by Corollary~\ref{corsx}}
\\ &\simeq& R_{(k)}[t]/(t^e-u) \otimes_{R_{(k)}} R_{(k)}/(u) \text{\ \ by Prop.~\ref{struc}}
\\& = & \bar\Q_p[t]/(t^e)
\end{eqnarray*}
On the other hand $\Symb_\Gamma^\pm(\D_k)_{(f_\beta)}= \Symb_\Gamma^\pm(\D_k)_{(x)}$ is free of rank one over $(\TT^\pm_{k,\nu})_{(x)}$ by Prop.~\ref{freerankone}. This proves the theorem.
\end{pf} 

\begin{cor} \label{corL} The $\HH$-eigenspace $\Symb_\Gamma^\pm(\D_k)[f_\beta]$ is $t^{e-1} \Symb_\Gamma^\pm(\D_k)_{(f_\beta)}$ and  
has dimension $1$ over $\bar \Q_p$. Its image by $\es^\pm \circ \rho_k^\ast$ is $0$ if $e>1$;
if $e=1$ it is the line $M_{k+2}(\Gamma)[f_\beta]$  when $f$ is cuspidal, or when $f$
is Eisenstein and the sign $\pm$ is $\epsilon(f_\beta)$, and $0$ when $f$ is Eisenstein 
and the sign $\pm$ is $-\epsilon(f_\beta)$.
\end{cor}
\begin{pf} Generally speaking, the eigenspace for some commuting linear operators and system of eigenvalues in a finite-dimensional space
is the {\bf socle} of the generalized eigenspace (for the same operators and eigenvalues) seen as a module over the algebra generated
by those operators. In our case, the eigenspace we are looking at is the socle of $\Symb_\Gamma^\pm(\D_k)$ as a module over
 $(\TT^\pm_{W,\rho})_{x}\simeq \bar \Q_p[t]/(t^e)$. This socle is one dimensional, generated by the image of $t^{e-1}$.

For the second part, we note that by Corollary~\ref{corcontrol}, the maps $\es^\pm \circ \rho_k^\ast$ induce  
Hecke-equivariant maps 
\begin{eqnarray*} \es^+ \circ \rho_k^{\ast} : \Symb_\Gamma^+(\D_k)_{(f_\beta)} &\rightarrow& M_{k+2}(\Gamma)_{(f_\beta)}\\
\es^- \circ \rho_k^{\ast} : \Symb_\Gamma^-(\D_k)_{(f_\beta)} &\rightarrow& M_{k+2}(\Gamma)_{(f_\beta)}
\end{eqnarray*} 
Note that the target of both maps are  one dimensional
since $f_\beta$ is a newform (see Lemma~\ref{lemmanewform}(i)), 
and that (by Corollary~\ref{corcontrol} and \S\ref{signeisenstein}) 
both maps are surjective except when $\pm=-\epsilon(f_\beta)$ and $f$ is Eisenstein, in which case the map is $0$.

 On the other hand the sources of both maps have dimension $e$ by the theorem above. If $e=1$, the  inclusion 
$$\Symb_\Gamma^\pm(\D_k)[(f_\beta)] \subset \Symb_\Gamma^{\pm}(\D_k)_{(f_\beta)}$$ 
 is an equality, so the maps $\es^\pm \circ \rho_k^\ast$ are surjective from $\Symb_\Gamma^\pm(\D_k)[(f_\beta)]$ to their respective targets (save the exceptional case, where the map is $0$). Therefore they are an isomorphism by equality of dimensions of the sources and targets, 
 except when the sign $\pm$ is $-\epsilon(f_\beta)$ and $f$ is Eisenstein, in which case the map is $0$.

In general, for $e \geq 1$, the operator $t$ acts nilpotently on the source of $\es^\pm \circ \rho_k^\ast$, so nilpotently on the target, 
but since this target has dimension at most $1$, it acts by $0$.
So if $e>1$, then $t^{e-1}$ acts by $0$ on $M_{k+2}(\Gamma)_{(f_\beta)}$, and the image of our eigenspace is $0$.
\end{pf}

We note that we have now proved Theorem~\ref{eigenspace}. Indeed, the first assertion of that theorem, namely the 
$1$-dimensionality of  $\Symb_\Gamma^\pm(\D_k)[f_\beta]$ is the first sentence of Corollary~\ref{corL}.
The second assertion, namely that the dimension of $\Symb_\Gamma^\pm(\D_k)_{(f_\beta)}$ is independent of the sign, and is equal to $1$
if and only if $f_\beta$ is non-critical follows from theorem~\ref{thmL}, which says that that dimension is the degree of ramification $e$ of
$\CC^0$ or $\CC$ over the weight space at $x$, and of Prop.~\ref{propcrit}(ii) which states that $e=1$
if and only if $f_\beta$ is non-critical. The last assertion of Theorem~\ref{eigenspace} follows immediately from the second sentence of Corollary~\ref{corL}.

\subsubsection{Properties of the $p$-adic $L$-function and proof of Theorem~\ref{thminterpol}}

\begin{definition} \label{defL} We now choose a non-zero vector $\Phi_{f_\beta}^\pm$ in $\Symb_\Gamma^\pm(\D_k)[f_\beta]$ (well defined up to a scalar by Corollary~\ref{corL}) and 
define the $p$-adic $L$-function
$$L^\pm(f_\beta,\sigma) = \Lambda_k (\Phi_{f_\beta}^\pm)(\sigma),$$
which clearly amounts to the same definition as in the introduction.
\end{definition}

\par \bigskip
We check the properties of those functions listed in the introduction: The {\bf analyticity} is clear, the estimation of the 
{\bf order} of $L^\pm(f_\beta,\cdot)$ follows from the result of Visik recalled in Prop~\ref{proporder}. 

\par \bigskip
The {\bf interpolation} properties follow immediately from Corollary~\ref{corL} and Prop.~\ref{inter0} in case $\rho_k^\ast(\Phi_{f_\beta}^\pm)=0$ and 
from Prop~\ref{interf} in case $\rho_k^\ast(\Phi_{f,\beta}^\pm) \neq 0$ and $f$ is 
cuspidal. 

\par \bigskip
The remaining case needed for the proof of Theorem~\ref{thminterpol}
is the case where $f$ is Eisenstein, and $\rho_k^\ast(\Phi_{f,\beta}^\pm) \neq 0$, which, by Corollary~\ref{corL}, implies that $\pm = \epsilon(f_\beta)$. This case needs a specific and careful treatment:

To ease notations, let us denote $\phi$ for $\Phi_{f_\beta}^{\epsilon(f_\beta)}$.
Let us call $\mu \in \D[0]$ the distribution attached to $\phi$ as in \S\ref{from},
that is $\tilde \phi (\{\infty\} -\{0\})$. Since $L^{\epsilon(f_\beta)}(f_\beta,.)$ is by definition the Mellin transform of the restriction of $\mu$ to $\Z_p^\ast$, the interpolation property would easily follow from the following statements:
\begin{eqnarray} \label{statint1} 
\forall j, \ 0 \leq j < k, \ \ \forall n \geq 1,\ \  \forall a=1,\dots,p^n-1,\  p \not | a,\, 
 \int_{\Z_p^\ast} \un_{a+p^n\Z_p}(z) z^j \ d\mu(z) = 0 \\ \label{statint2}
  \forall n \geq 1,\ \  
 \int_{\Z_p^\ast} \un_{a+p^n\Z_p}(z) z^k \ d\mu(z) \text{\ is independent of }a=1,\dots,p^n-1, p\not |a, 
\end{eqnarray}

Since $\rho_k^\ast(\phi)$ is an Eisenstein symbol, there exists by Prop~\ref{es} 
an element
$\tau \in \Symb_\Gamma(\Delta)$ ($\Delta$ the abelian groups of divisors on $\PP^1(\Q)$, not necessarily of degree $0$), such that 
\begin{eqnarray} \label{rhotau} \rho^\ast(\phi)(\{\infty\}-\{a/p^n\})=\tau(\{\infty\})-\tau(\{a/p^n\}) \end{eqnarray} 
Let us write $\tau(\{\infty\})=\sum_{i=0}^{k} d_i X^i Y^{k-i}$ and $\tau({a/p^n})=
\sum_{i=0}^d e_i X^i Y^{k-i}$ (note that $e_i$ depends on $a \in \{1,\dots,p^n-1\}$), so that we have, by equation (\ref{calculbis})
of Lemma~\ref{interpolation} 
\begin{eqnarray} \label{di} \int z^j {\bf{1}}_{a+p^n\Z_p} \ d \mu(z) & = &
\beta^{-n} \sum_{i=0}^j d_{i} (-1)^{k-i} \frac{j!(k-i)!}{k!(j-i)!} p^{ni}a^{j-i} \\
\label{ei}
& - &  \beta^{-n} \sum_{i=0}^j e_{i} (-1)^{k-i} \frac{j!(k-i)!}{k!(j-i)!} p^{ni}a^{j-i}
\end{eqnarray}
We study both terms (\ref{di}) and (\ref{ei}) separately, showing that they both can be non-zero only when $j=k$, and that they are at any rate independent on $a$. 

For (\ref{di}), the fact that $\tau(\{\infty\})$ is invariant by the stabilizer
of $\{\infty\{$ in $\Gamma$ shows that this polynomial is of the form $d_k X^k$, with all $d_j =0$ for $j <k$. Therefore, the term (\ref{di}) is non $0$ only when $j=k$. It is independent on $a$ since the $d_i$ are and since all the terms with $i < j=k$ in the sum are $0$.

 As for (\ref{ei}), the fact that $\tau(\{a/p^n\})$ is invariant by the stabilizer
of $\{a/p^n\}$ in $\Gamma$ shows that this polynomial has the form $e (aX+p^nY)^k$, with $e$ some constant, possibly depending on $a$. If $k>0$, a simple computation shows that the term (\ref{ei}) is always $0$ (for all values of $j$, $0 \leq j \leq k$), which is what we wanted. If $k=0$, the value of $(\ref{ei})$ is $e= \tau(\{a/p^n\})$ and we have to prove that this value is independent on $a=1,\dots,p^n-1$ for $a$ prime to $p$. 
Since $N$ is prime to $p$, all
the cusps $\{a/p^n\}$ for $a$ as above are in the same 
$\Gamma_1(N) \cap \Gamma_0(p)$-orbits 
(cf \cite[Proof of prop. 1.43]{shimura}, or \cite[Proof of Theorem 4.2.9]{miyake}),
sau $\tau(\{a/p^n\})$ is independent of $a$. 

We thus have proved (\ref{statint1}) and (\ref{statint2}),  hence Theorem~\ref{thminterpol}.

\begin{example}\label{exLEisenstein} Take $N=1$ and $k > 0$ even. As is well known, there exists one normalized Eisenstein series $f=E_{k+2,1,1}$ of weight $k+2$ and level $1$. Let $f_\beta$ (with $\beta=p^{k+1}$) be its critical-slope refinement. The only non-zero boundary eigensymbol  $\tau \in \Hom(\Delta,V_k)^\Gamma$ for $\Gamma=\Gamma_0(p)$ and of weight $k+2$ satisfying $U_p \tau = p^{k+1} \tau$ is, up to a scalar, the one such that $\tau(\infty)=X^k$, $\tau(0)=\frac{p^{k+1}-1}{p^k(p-1)} Y^k$. By Prop~\ref{es}, we have $\tau=\rho_k^\ast(\Phi^{\epsilon(f_\beta)}_{f_\beta}$, and it easy to see that actually $\epsilon(f_\beta)=1$. Then the above computations show that, with $\mu$ as above, $\int z^j {\bf{1}}_{a+p^n\Z_p} \ d \mu(z) = \frac{\delta_{j,k}}{p^n}$. Therefore, $\int z^k d\mu(z) = (p-1)/p$  is non $0$, so $L^+(f_\beta,t^k) \neq 0$.
\end{example}
\par \bigskip

Finally, it remains to prove Corollary~\ref{corinterpol},  that is that for $f_\beta$ of critical slope, $L^\pm(f_\beta,\sigma)$ has infinitely many $0$'s in every open unit disc (that is, connected component) of $\WW$. When $f_\beta$ is $\theta$-critical, the interpolation properties already gives an infinity of $0$. When it is not, the proof follows exactly the one given in \cite{pollack} (cf. Theorem 3.3. there, attributed to Mazur): For any fixed  character $\psi$ of $\mu_{p-1}$
we observe from the interpolation properties the existence
of elements $c_n \in \C_p^\ast$ such that for $n \geq 1$
\begin{eqnarray*} 
L_p(f_\beta, \psi \chi_{\zeta_n-1})& =& \frac{c_n}{\beta^{n+1}} \\
L_p(f_\alpha, \psi \chi_{\zeta_n-1}) &=& \frac{c_n}{\alpha^{n+1}}
\end{eqnarray*}
where $\alpha$ is the other root of (\ref{hp}), $\zeta_n$ is a primitive $p^n$-root of $1$ in $\bar \Q_p$,
and as in the introduction $\chi_{\zeta_n-1}$ is the character of $1+p\Z_p$ that sends the chosen generator $\gamma$ of the group to $\zeta_n-1$.
As shown in \cite[Theorem 3.3]{pollack}, this implies that on each connected component of $\WW$, either $L(f_\beta,.)$ or $L(f_\alpha,.)$ has an infinity of $0$'s. But $v_p(\alpha)=0$ since $v_p(\beta)=k+1$, so $L(f_\alpha,.)$ has order $0$, and
thus has only a finite number of $0$'s. This completes the proof.

\subsection{Two-variables $p$-adic $L$-functions}

We keep the notations of the preceding subsection.

\subsubsection{Modular symbols over the eigencurve}

\begin{prop} \label{struc2} 
Up to shrinking the affinoid neighborhood $W=\sp R$ of $k$ in $W$, there exists an affinoid neighborhood $V =\sp T$ of 
$x \in \CC^\pm_{W,\nu}$ which is a connected component of $\sp \TT_{W,\nu}$, an element $u \in R \otimes \bar \Q_p$ such that $u(k)\neq 0$, and $k$ is the only $0$ of $u$ on $U$, 
an element $t \in T$ such that $t(x) \neq 0$,
and an isomorphism of $R$-algebras $T \simeq R[X]/(X^e-u)$ sending $t$ on $X$.
\end{prop}
\begin{pf} 
This follows easily from Prop~\ref{struc}.
\end{pf}

We consider the $T$-module $$M := \Symb_\Gamma^\pm(\D \otimes R)^{\leq \nu} \otimes_{\TT^\pm_{W,\nu}} T.$$  
It is free of rank one.
Note that $T = \epsilon \TT^\pm_{W,\nu}$ for some idempotent $\epsilon$, so $M:= \epsilon  \Symb_\Gamma^\pm(\D \otimes R)$
is a submodule of  $\Symb_\Gamma^\pm(\D \otimes R)$.

Now we set $$N = M \otimes_R T.$$ This $R$-modules has two structures of $T$-module, one coming from the $T$-structure that $M$ already
had (called below {\it  the first $T$-structure}), the other coming from the right factor $T$ in the tensor product (called {\it the second 
$T$-structure}). 

We might think of $N$ as the module of distributions-valued modular symbols over the eigencurve $\CC$ (more precisely, over the open subspace $V$ of $\CC$)
as opposed as $M$ which is (a part of) the module of distributions-valued modular symbols  over the weight space (more precisely, over the open subset $W$ of $\WW$).

\par \medskip
For every point $y \in V(\bar\Q_p)$, of weight $\kappa(y) \in W(\bar \Q_p)$ we define 
a {\it specialization map} 
$$\spe_y :\, N \rightarrow \Symb_\Gamma^\pm(\D_{\kappa(y)}) \otimes_{\Q_p} {\bar \Q_p}$$
by composing the following canonical morphisms
\begin{eqnarray*} N & \rightarrow & N \otimes_{T,y} \bar \Q_p \text{\ \ \ where the tensor product is for the second $T$-structure}\\ 
 &=& (M \otimes_R T) \otimes_{T,y} \bar \Q_p  \\
&=& M \otimes_{R,\kappa(y)} \bar \Q_p  \\
& \subset & \Symb_\Gamma^\pm(\D_{\kappa(y)}) \otimes_{\Q_p} {\bar \Q_p}
\end{eqnarray*}
This map is $\HH$-equivariant if we give $M$ the {\it first} $T$-structure (and the corresponding action of $\HH$) and 
$\Symb_\Gamma^\pm(\D_{\kappa(y)})$ its natural action of $\HH$.

\subsubsection{Mellin transform over the eigencurve}

Recall from~\ref{mellin} that there is an $R$-linear map 
$$\Lambda_R : \Symb_\Gamma(R \hotimes \D)^{\leq \nu} \rightarrow R \hotimes \RR.$$
By restriction this maps induces a map
$$\Lambda_{|M} :M \rightarrow R \hotimes \RR$$
from which we get a map
$$\Lambda_T:=\Lambda_{|M} \otimes \Id_T : N=M \otimes_R T \rightarrow (R \hotimes \RR) \otimes_R T \isomo T \hotimes \RR.$$
Note that elements of the target of this map are two-variables analytic functions, the first variable being a point 
$y \in \sp T = V \subset \CC$, the second variable being an element $\sigma$ of $\WW$.

The morphism $\Lambda_T$ is obviously $T$-equivariant if we give $N$ its second $T$-structure.

\begin{lemma} \label{compLambdaT} If $\Phi \in N$ and $y \in V(\bar \Q_p)$, then for all $\sigma$, we have 
$$\Lambda_T(\Phi)(y,\sigma)=\Lambda(\spe_y(\Phi))(\sigma).$$
\end{lemma} 
\begin{pf}
Assume first that $\Phi$ is a pure tensor $\phi \otimes t$ with $\phi \in M$ and $t \in T$.
Then $\spe_y(\Phi)=\spe_y(\phi \otimes t) = t(y) \phi \in M \subset \Symb_\Gamma^{\pm} (\D_{\kappa(y)})$ by definition of $\spe_y$.
Therefore $\Lambda(\spe_y(\Phi))(\sigma) = t(y) \Lambda(\phi)(\sigma)$ for all $\sigma \in \WW(\C_p)$.
On the other hand $\Lambda_T(\Phi)=\Lambda_T(\phi \otimes t) = \Lambda_R(\Phi) \otimes t \in (R \hotimes \RR) \otimes T \isomo T \hotimes R$.  Therefore $\Lambda_T(\phi)(y,\sigma) = (\Lambda_R(\phi) \otimes t)(y,\sigma) = t(y) \Lambda(\phi)(\kappa(y),\sigma)$. 
The equality $\Lambda(\spe_y(\Phi))(\sigma)=\Lambda_T(\Phi)(y,\sigma)$ then follows from (\ref{complambda}) applied with $w=\kappa(y)$.

The case where $\Phi$ is a general element of $N$, that is an $R$-linear combination of pure tensors as above, follows immediately if we observe that the two maps $N \rightarrow \RR$ given by $\Phi \mapsto \Lambda_T(\Phi)(y,\cdot)$ and $\Phi \mapsto \Lambda(\spe_y(\Phi))$ are both $R$-linear if we give $N$ its natural $R$ structure (this is the same for the two $T$-structure) and
$\RR$ its $R$-structure using the morphism $\kappa(y): R \rightarrow \bar \Q_p$.
\end{pf}

\subsubsection{Proof of Theorem~\ref{twovariables}}

\label{proof2} 
Let us choose $\phi$ a generator of the free $T$-module $M$, and set $$\Phi = \sum_{i=0}^{e-1} t^i \phi \otimes t^{e-1-i} \in N.$$
It is clear that the element $\Phi$ of $N$ depends on the choices made (the isomorphism in Prop~\ref{struc2}, the generator $\phi$ of $M$) only up to multiplication by an element of $T^\ast$ for the first $T$-structure on $N$.

\begin{lemma} \label{lemmaPhi} One has $$(t \otimes 1 - 1 \otimes t) \Phi = 0.$$ In particular, the action of any element of $T$ on $\Phi$ for the first or the second $T$-structure on $N$ are identical.
\end{lemma}
\begin{pf}
We compute in $N=M \otimes_R T$,
\begin{eqnarray*} 
(t \otimes 1 - 1 \otimes t) \Phi &=&  (t \otimes 1 - 1 \otimes t)  \sum_{i=0}^{e-1} t^i \phi \otimes t^{e-1-i} \\
 & = & \sum_{i=0}^{e-1} t^{i+1} \phi \otimes t^{e-1-i} - t^i \phi \otimes t^{e-i} \\
 & = & t^e \phi \otimes 1 - \phi \otimes t^e \text{ \ \ \ "telescopic cancellations"}\\
 &=& t^e(\phi \otimes 1 - \phi \otimes 1) \text{\ \ \ since $t^e = u$ is in $R$}\\
 &=& 0
 \end{eqnarray*}
\end{pf}
In particular, the element $\Phi$ is well-defined up to a multiplication by an element of $T^\ast$ for the second $T$-structure 
as well.

\begin{prop} \label{prop2var}
If $y \in V(\bar \Q_p)$ corresponds to a refined modular form $f'_{\beta'}$ of weight $k'=\kappa(y)$, 
then in $\Symb_\Gamma^\pm(\D_{k'}) \otimes_{\Q_p} \bar \Q_p$
we have, up to a scalar in $\bar \Q_p^\ast$, the equality
$$\spe_y (\Phi) = \Phi_{f'_{\beta'}},$$
where $\Phi_{f_{\beta'}}$ is any generator of the one dimensional $\bar \Q_p$-space $\Symb_\Gamma^\pm(\D_{k'})_{(f'_{\beta'})}$
(see Corollary~\ref{corL}) 
\end{prop}
\begin{pf}
It is enough, by Corollary~\ref{corL} to prove that $\spe_y(\Phi)$ is in the $\HH$-eigenspace in $M \otimes_{R,k'} \bar \Q_p \subset \Symb_\Gamma^\pm(\D_{k'})$
for the same system of $\HH$-eigenvalues as $f'_{\beta'}$. But the algebra of Hecke operators acting on
$M \otimes_{R,k'} \bar \Q_p$ is generated by a single operator $t$, and we thus only have to show that $t$
acts on $\spe_y \Phi$ as in acts on $f'_{\beta'}$, that is with the eigenvalues $t(y)$. This is clear by Lemma~\ref{lemmaPhi} and the 
definition of $\spe_y$.
\end{pf}

\begin{definition} \label{defLL}
We set $L^\pm := \Lambda_T(\Phi) \in T \otimes \TT$ and call it the {\it two variables $p$-adic $L$-function}.
It is an analytic function of two variables $(y,\sigma) \mapsto L^\pm(y,\sigma)$
defined on $V \otimes \WW$. The function $L^\pm$ is well-defined up to multiplication by a non-zero analytic function on $V$ (that is an element of $T^\ast$).
\end{definition}

If now $y$ is a point of $V(\bar \Q_p)$ corresponding to a refined modular form $f'_{\beta'}$ of weight $k'$, then 
\begin{eqnarray*}
L^\pm(y,\sigma) & = & \Lambda_T(\Phi)(y,\sigma) \text{\ \  by Definition~\ref{defLL}} \\
&=& \Lambda( \spe_y(\Phi)(\sigma) ) \text{\ \ by Lemma~\ref{compLambdaT}} \\
&=& \Lambda( \Phi_{f'_{\beta'}})(\sigma )\text{ \ \ (up to a scalar) \ \ by Prop.~\ref{prop2var}}\\
&=& L^\pm(f'_{\beta'},\sigma) \text{\ \ (up to a scalar) \ \ by Definition~\ref{defL}}
\end{eqnarray*}

This completes the proof of Theorem~\ref{twovariables}

\subsection{Secondary $p$-adic $L$-functions}

 We keep all notations and assumptions as above. We assume moreover than $e>1$, that is that $f_\beta$ is critical.
 Since $t$ is a uniformizer of $V$ at $x$, we can consider the two variable $L$-function $L(y,\sigma)$ with $y$ in a neighborhood of $x$ as $L(t,\sigma)$ with $t$ in a neighborhood of $0$. 
 
\begin{definition} For $i=0,\dots,e-1$,
 $$L_i^\pm(f_\beta,\sigma) =  \left(\frac{\partial^i L^\pm(t,\sigma)}{\partial{t^i}}\right)_{|t=0}.$$
 \end{definition}
 This is in accordance with the definition given in the introduction.
 
 To get more information on those secondary $L$-functions, we compute 
 \begin{eqnarray*}
 L^\pm(t,\sigma) &=& (\Lambda_T) (\sum_{i=0}^{e-1} t^i \phi \otimes t^{e-1-i}) \text{ \ \ \ by Definition~\ref{defLL}}
\\
&=& \sum_{i=0}^{e-1}  t^{e-1-i} \Lambda_T (t^i \phi \otimes 1) \text{\ \ \ since $\Lambda_T$ is $T$-linear for the second $T$-structure on $N$} 
\end{eqnarray*}
from which we get:
\begin{eqnarray*} L_i^\pm(f_\beta, \cdot) &=& \left(\frac{\partial^i L(t,\cdot)}{\partial{t^i}}\right)_{t=0}  \\
& = &  i! \Lambda_T(t^{e-1-i}\phi \otimes 1)_{|t=0} \text{\ \  by Leibniz' rule} \\
&=& i! \Lambda( \spe_x ( t^{e-1-i} \phi \otimes 1) ) \text{\ \ by Lemma~\ref{compLambdaT}} \\
&=&  i! (\Lambda t^{e-1-i} \phi_{k} )
\end{eqnarray*}
 where $\phi_{k}$ is the image of $\phi$ in $M \otimes_{R,k} \bar \Q_p = \Symb_\Gamma^\pm(\D_k)_{(x)}$.
 Since $\phi$ is an $T$-generator of $M$, $\phi_{k}$ is also a $T \otimes_{R,k} \bar \Q_p = (\TT_{k,\nu}^\pm)_{(x)}$-generator
 of  $\Symb_\Gamma^\pm(\D_k)_{(x)}$. Recall that by theorem~\ref{thmL}, $(\TT_{k,\nu}^\pm)_{(x)}$ is isomorphic to $\bar \Q_p[t]/t^e$
 and the module $\Symb_\Gamma^\pm(\D_k)_{(x)}$ is free of rank one over that algebra. So the flag 
 $(t^{e-1} \phi_k, t^{e-2} \phi_k, \dots, \phi_k)$ is independent of the generator $\phi_k$ in $\Symb_\Gamma^\pm(\D_k)_{(x)}$
 hence independent of $\phi_k$. As we have seen, the image by $\Lambda$ of the generalized eigensymbols 
 $t^{e-1} \phi_k, t^{e-2} \phi_k, \dots, \phi_k$ are up to a scalar the secondary $L$-function $L^\pm(f_\beta,\cdot), L_1^\pm(f_\beta,\cdot),\dots, L_{e-1}^\pm (f_\beta, \cdot)$ we have proved the {\bf connection with modular symbols} of the  secondary $L$-functions
 stated in the introduction. 
 
 We now prove the other properties of the $L_i^\pm(f_\beta,\cdot)$ The {\bf Order} properties follows directly from~\ref{proporder}.
 
 The {\bf interpolation} properties for $L_i^\pm(f_\beta,\cdot)$ for $0 \leq i \leq e-2$ follows from Prop.~\ref{inter0}, since
 that function is up to a scalar $ \Lambda(t^{e-1-i} \phi_k) $ and $\rho_k^\ast(t^{e-1-i} \phi_k) =  t^{e-1-i} \rho_k^\ast (\phi_k)= 0$ for $i < e-1$ since $t$ acts by $0$ in $\Symb_\Gamma^\pm(V_k)_{(x)}$ (see the proof of Corollary \ref{corL}). However, when $i=e-1$, one finds that
$L_{e-1}^\pm(f_\beta,\cdot) = \Lambda( \phi_k)$ and $\phi_k$ is a generator of $\Symb_\Gamma^\pm(V_k)_{(x)}$, so since $\rho_k^\ast$ is surjective with non $0$ image,
we have $\rho_k^\ast(\phi_k) \neq 0$. The {\bf interpolation} property thus follows from Prop.~\ref{interf}. 
 
 The {\bf infinity of 0's} properties for $L_i(f_\beta,\cdot)$ for $i < e-1$ is a trivial consequence of the interpolation property,
 while for $L_{e-1}(f_\beta,\cdot)$ it follows exactly as in the case of a non-critical form, considering the pair of functions
 $L_{e-1}(f_\beta,\cdot)$ and $L(f_\alpha,\cdot)$.

\par \bigskip
\end{document}